\documentclass[10pt]{amsart}

\usepackage{amsmath,amssymb,latexsym,soul,cite,amsthm,color,enumitem,graphicx,tikz, mathtools,microtype}
\usepackage[utf8]{inputenc}
\usepackage{amsfonts,amsmath,amsthm,latexsym,amssymb,marvosym,esint,xfrac,bbm,amssymb}
\usepackage{xcolor}
\usepackage{tikz}
\tikzstyle{mybox} = [draw=black, very thick, rectangle, rounded corners, inner ysep=5pt, inner xsep=5pt]
\usepackage{hyperref}
\usepackage{verbatim}
\usepackage{ulem}
\usepackage[makeroom]{cancel}
\usepackage[margin=1in]{geometry}
 \usepackage[usenames,dvipsnames]{pstricks}
 \usepackage{pst-grad} 
 \usepackage{pst-plot} 
\allowdisplaybreaks

\DeclareMathOperator*{\osc}{osc}

\numberwithin{equation}{section}
\newtheorem{theorem}{Theorem}[section]
\newtheorem{lemma}{Lemma}[section]
\newtheorem{remark}{Remark}[section]

\newtheorem{propo}{Proposition}[section]

\def\XXint#1#2#3{{\setbox0=\hbox{$#1{#2#3}{\int}$}
     \vcenter{\hbox{$#2#3$}}\kern-.5\wd0}}

\title{Qualitative properties of solutions to anisotropic parabolic equations: Part III- H\"{o}lder continuity
}
\author{ Simone Ciani, Eurica Henriques, Maria O. Savchenko, Igor I. Skrypnik
 }
\begin{document}

\begin{abstract}

\noindent In this work, we study the local H\"older continuity of non-negative solutions to a class of doubly nonlinear anisotropic partial differential equations, with prototype
\[
u_t-\sum_{i=1}^N D_i\Big(|D_i u^{m_i}|^{p_i-2}D_i u^{m_i}\Big)=0,
\quad m_i>0,\quad p_i>1,\quad i=1,\ldots,N. 
\]
The main result shows that, when the exponents describing the diffusion in different spatial directions are sufficiently close to one another, the solutions enjoy local H\"older continuity. A key feature of the analysis is the distinction between two diffusion regimes: the degenerate and the singular cases; the local regularity results are established in both regimes, under suitable conditions on the exponents.

\vspace{.2cm}


\end{abstract}

\maketitle

\begin{center}
		\begin{minipage}{9cm}
			\small
  \tableofcontents
		 \end{minipage}
	 \end{center}

\addtocontents{toc}{\protect\setcounter{tocdepth}{1}}
    
\section{Introduction and Main Results}\label{Introduction}

\noindent In this paper we extend the classical local regularity theory to a class of doubly nonlinear anisotropic parabolic equations whose prototype is of the form
\[ u_t-\sum\limits_{i=1}^N D_i \Big(u^{(m_i-1)(p_i-1)} |D_i u|^{p_i-2} D_i u \Big)=0,\quad m_i>0,\quad p_i>1,\,\,i=1,..., N,\,\,u\geqslant 0. \]
This (now) classical regularity theory for nonlinear differential equations dates back  to the works DiBenedetto  \cite{DiB_1, DiB_2} and Chen and DiBenedetto \cite{CheDiB}, where they established the interior H\"{o}lder continuity of solutions to quasi-linear parabolic equations of $p$-Laplace type with merely measurable coefficients; as for the porous medium type equations, analogous regularity results were obtained by Caffarelli and Friedman \cite{CafFri}, Caffarelli and Evans \cite{CafEva}, and DiBenedetto and Friedman \cite{DiBFri}. These foundational contributions were later generalized to the framework of doubly nonlinear parabolic equations, beginning with the pioneering works of Ivanov \cite{Iva1, Iva2}, Porzio and Vespri \cite{PorVes}, and Vespri \cite{Ves1, Ves2}. 

\noindent When considering elliptic and parabolic anisotropic contexts, the boundedness of solutions has been well understood since the seminal works of Kolodii \cite{Kol1, Kol2} (see also \cite{CiaHenSavSkr1} and \cite{CVV} for recent results on this topic), however finer qualitative properties - such as continuity and Harnack-type inequalities - remain far less explored in the anisotropic setting. Concerning $C^0$-regularity, the anisotropic porous medium equation has attracted the most attention; see, for instance, \cite{BurSkr, Hen1, Hen2}. To the best of our knowledge, the available results for anisotropic elliptic equations of the form
\[
\sum\limits_{i=1}^N D_i \Big(|D_i u|^{p_i-2} D_i u \Big)=0,\quad p_i>1,\,\,i=1,..., N,
\] 
are still partial and rely on restrictive structural assumptions. Specifically, results in \cite{LiaSkrVes, CiaSkrVes, LisSkr, CiaHenSkr} are obtained under conditions such as $p_1=2 < p=p_2=\cdots=p_N$ or $p=p_1=\cdots=p_s<2=p_{s+1}= \cdots=p_N$ or $p=N$, 
whereas in \cite{DiBGiaVes4} the analysis requires the additional smallness condition $p_+-p_- \ll 1$. In the parabolic anisotropic setting, only partial results are currently available, notably Harnack-type inequalities in certain model cases such as the ones presented in \cite{CiaGua, CiaMosVes, CiaHen, CiaHenSkr24}.

\noindent The present work continues the investigation of regularity properties for anisotropic parabolic equations initiated in \cite{CiaHenSavSkr1, CiaHenSavSkr2, CiaHenSavSkr3}. Our main contribution is the establishment of local H\"{o}lder continuity for the non-negative, locally bounded, local weak solutions under general structural conditions. More precisely, we will consider a non-negative function $u$, defined in a cylindrical domain $\Omega_{T}:= \Omega \times (0, T]$, where $\Omega$ is a bounded domain in $\mathbb{R}^{N}$ and $T>0$, satisfying 
\begin{equation}\label{eq1.1}
u_{t}-\sum\limits_{i=1}^N D_i a_i (x, t, u, D u)=0, \quad (x, t)\in \Omega_{T}
\end{equation}
being $a_i:\Omega_{T}\times \mathbb{R}_+\times \mathbb{R}^{N} \rightarrow \mathbb{R}^{N}$ such that $a_i(\cdot, \cdot, u, \xi)$ are Lebesgue measurable for all $u\in \mathbb{R}_+$, $\xi \in \mathbb{R}^{N}$, while $a_i(x, t, \cdot, \cdot)$ are continuous for almost every $(x, t)\in \Omega_{T}$, and verify the following structure conditions:
\begin{equation}\label{eq1.2}
\begin{cases}
\sum\limits_{i=1}^N a_i(x, t, u, D u)\,D_i u^{m_-}  \geqslant  K_{1}\,\sum\limits_{i=1}^N u^{(m_i-m_-)(p_i-1)}\,|D_i u^{m_-}|^{p_i},\\
\displaystyle{|a_i(x, t, u,  D u)| \leqslant  K_{2}\,  u^{\frac{(m_i-m_-)(p_i-1)}{p_i}} \left(\sum_{j=1}^N u^{(m_j-m_-)(p_j-1)}|D_j u^{m_-}|^{p_j}\right)^{\frac{p_i-1}{p_i}}, \quad  i=1, ...,N,}
\end{cases}
\end{equation}
where $K_1$, $K_2$ are positive constants, $p_i>1$, $m_i>0$ and $m_-= \min\limits_{1\leqslant i\leqslant N} m_i$.

\noindent The ongoing study requires a careful understanding of the values of the various exponents at hands, the adopted geometry, the functional setting framework and the meaning of local weak solution. That is precisely the content of the next lines. 

\noindent {\bf Exponents.} For fixed numbers $p_1,..., p_N >1$, $m_1, ..., m_N >0$ we set
$$p_+:=\max\limits_{1\leqslant i\leqslant N}p_i,\qquad p_-:=\min\limits_{1\leqslant i\leqslant N}p_i,\qquad \frac{1}{p}:=\frac{1}{N}\sum\limits_{i=1}^N\frac{1}{p_i}, \ \  p<N $$
$$m_+:=\max\limits_{1\leqslant i\leqslant N}m_i,\qquad m_-:=\min\limits_{1\leqslant i\leqslant N} m_i \ ,$$
$$\lambda_i:=m_i(p_i-1),\qquad \lambda_+:=\max\limits_{1\leqslant i\leqslant N} \lambda_i,\qquad \lambda_-:=\min\limits_{1\leqslant i\leqslant N}\lambda_i\ ,$$
$$ \Big|\frac{\lambda}{p}\Big|:=\frac{1}{N}\sum\limits_{i=1}^N \frac{\lambda_i}{p_i} , \qquad \varkappa = p+N \left( p\Big|\frac{\lambda}{p}\Big| -1\right)  \ .$$

\noindent {\bf Geometric setting.} For  $\vec{r}:=(r_1, ..., r_N)$, $\eta >0$ we define the 
\begin{enumerate}
    \item[] anisotropic cube
\[K_{\vec{r}}(\bar{x}):=\big\{x:|x_i-\bar{x}_i|<r_i,\quad i=1, ...,N\big\} \ , \]
and the corresponding anisotropic cylinder
\[ Q_{\vec{r}, \eta}(\bar{x}, \bar{t}):=K_{\vec{r}}(\bar{x})\times (\bar{t}-\eta, \bar{t}) \] 
        \item[] intrinsic anisotropic cube
\[ K^k_r(\bar{x}):=\big\{x: |x_i-\bar{x}_i|<r^{\frac{p_+}{p_i}} k^{-\frac{\lambda_+-\lambda_i}{p_i}},\quad i=1, ..., N\big\},\quad k,\,r>0 \ , \]
and the intrinsic anisotropic cylinder
\[ Q^k_{r, \eta}(\bar{x}, \bar{t}):=K^k_r(\bar{x})\times (\bar{t}-\eta, \bar{t}),\quad k,\,r,\,\eta>0 \ . \]
\end{enumerate}

\noindent {\bf Functional anisotropic framework.} We define the anisotropic spaces of locally integrable functions as
$$W^{1,\textbf{p}}_{loc}(\Omega):=\big\{u\in W^{1,1}_{loc}(\Omega): D_i u\in L^{p_i}_{loc}(\Omega),\,\,i=1, ..., N\big\},\quad \textbf{p}=(p_1, ..., p_N),$$
$$L^{\textbf{p}}_{loc}(0, T; W^{1, \textbf{p}}_{loc}(\Omega)):=\big\{u\in L^1_{loc}(0, T; W^{1,1}_{loc}(\Omega): D_i u\in L^{p_i}_{loc}(0, T; L^{p_i}_{loc}(\Omega)),\,\,i=1, ..., N\big\},$$
and the respective spaces of functions with zero boundary data
$$W^{1,\textbf{p}}_{0}(\Omega):=\big\{u\in W^{1,1}_{0}(\Omega): D_i u\in L^{p_i}_{loc}(\Omega),\,\,i=1, ..., N\big\},\quad \textbf{p}=(p_1, ..., p_N),$$
$$L^{\textbf{p}}_{loc}(0, T; W^{1, \textbf{p}}_{0}(\Omega)):=\big\{u\in L^1_{loc}(0, T; W^{1,1}_{0}(\Omega): D_i u\in L^{p_i}_{loc}(0, T; L^{ p_i}_{loc}(\Omega)),\,\,i=1, ..., N\big\}.$$

\noindent {\bf Local weak solution.} A function
\[
\begin{cases}
u\in C_{loc}(0, T; L^{m_-+1}_{loc}(\Omega)),\,\, u^{m_-}\in L^{\textbf{p}}_{loc}(0, T; W^{1,\textbf{p}}_{loc}(\Omega)),\\
u^{\frac{\lambda_i+m_-}{p_i}},\, D_i u^{\frac{\lambda_i+m_-}{p_i}} \in L^{p_i}_{loc}(0, T; W^{1, p_i}_{loc}(\Omega)),
\end{cases}
\]
is said to be a non-negative local weak sub (super)-solution to \eqref{eq1.1} if, for every compact set $E\subset \Omega$ and every subinterval $[t_1, t_2]\subset (0, T]$, the inequality
\begin{equation}\label{eq1.3}
\int\limits_E u\,\varphi\,dx \Big|^{t_2}_{t_1}+\int\limits^{t_2}_{t_1}\int\limits_E \Big[ - u\, \varphi_t+\sum\limits_{i=1}^N a_i(x, t, u, D u)\,D_i \varphi \Big]\,dx dt \leqslant (\geqslant) 0
\end{equation}
holds for all non-negative test functions
\[
\varphi\in W^{1,1+\frac{1}{m_-}}_{loc}(0, T; L^{1+\frac{1}{m_-}}(E))\cap L^{\textbf{p}}_{loc}(0,T; W^{1,\textbf{p}}_0(E)).
\]
A function is called a local weak solution if it is both a sub- and a super-solution.

\noindent {\bf Main Results.} As mentioned earlier, our goal is to present and prove results concerning the local H\"older continuity of the non-negative, locally bounded, local weak solutions \eqref{eq1.1}-\eqref{eq1.2}: we distinguish between the singular regime, $\lambda_+<1$, and the degenerate one,  $\lambda_+>1$. 

\begin{theorem}\label{th1.1}
Let $u$ be a non-negative, locally bounded, local weak solution to \eqref{eq1.1}-\eqref{eq1.2}. Consider  $\lambda_+<1$, $\lambda_+-\lambda_- \leqslant \epsilon_*$, for some sufficiently small $\epsilon_*$, and $\varkappa>0$. In this setting, $u$ is locally H\"{o}lder continuous in $\Omega_T$.
\end{theorem}

\begin{theorem}\label{th1.2} Let $u$ be a non-negative, locally bounded, local weak solution to \eqref{eq1.1}-\eqref{eq1.2} and consider  $\lambda_+>1$, $\lambda_+-\lambda_- \leqslant \epsilon_*$, $p_+-p_- \leqslant \epsilon_*^{(1)}$ for some sufficiently small $\epsilon_*, \epsilon_*^{(1)}$,  and $\bar{\varkappa}:=p +N(p-2)>0$, then $u$ is locally H\"{o}lder continuous in $\Omega_T$. 
\end{theorem}

\noindent The proofs of both results rely on deriving iterative schemes that lead to the reduction of the oscillation of u in a decreasing sequence of nested cylinders. The predominant tools involve expansion of positivity results together with DeGiorgi type lemmas (collected in the Subsections \ref{SS:resultsavailable} and \ref{SS:DEGlemmas}) and some integral inequalities of Harnack type which will be stated and proved along Section \ref{S:harnackineq}. One interesting feature is that the proof of Theorem \ref{th1.1}, presented in Section \ref{S:provingthe1.1}, is that it is, in a certain sense, simpler - one just need to consider the singular range $\lambda_+<1$;  the extra restriction on the exponents $\varkappa: =p+N\big(p\big|\frac{\lambda}{p}\big|-1\big)>0$ appears when  bringing measure theoretical information up in time. The proof of Theorem \ref{th1.2} (see Section \ref{S:provingthe1.2}), when considering the degenerate regime $\lambda_+>1$, is more demanding: one needs to distinguish between the degenerate regime regarding the $p_i$, $p_+>2$, and the singular one, $p_+<2$. In the late case, the additional restriction $\bar{\varkappa}:=p +N(p-2)>0$ is to be considered in order to make use of an integral Harnack type estimate. The results at hands require the exponents $\lambda_i$, $p_i$, $m_i$ not to be too sparse from each other.

\noindent {\bf Anisotropic (p)-Laplacian.} In the particular case \(m_i=1\), for all \(i=1,\ldots,N\), we recover the general evolutionary anisotropic (p)-Laplacian equation \eqref{eq1.1}, where the principal part of the differential operator satisfies
\[ a_i(x, t, u, D u)\,D_i u  \geqslant  K_{1}\, |D_i u|^{p_i},\quad 
|a_i(x, t, u,  D u)| \leqslant  K_{2}\,\sum_{i=1}^N |D_iu|^{\frac{p_i-1}{p_i}}, \ \  i=1, ...,N . \]
 The H\"older continuity has been established through different geometric and technical approaches in the degenerate (\(p_+>2\)) and singular (\(p_+<2\)) regimes, in line with the corresponding results available for the isotropic case. Our regularity results require the exponents to be sufficiently close to one another and rely on the condition \(p+N(p-2)>0\)
which ensures that an integral Harnack-type estimate can be established.

\section{Auxiliary Results}
\subsection{Known results available in the literature}\label{SS:resultsavailable}

\noindent In what follows we collect several known results which will be used along the text, just to keep this work as self contained as possible. We start by presenting two algebraic results that can be found in \cite[Lemma 3.1]{BogDuzGiaLiaSch} and \cite[Lemma 3.5]{CVV}, respectively.

\begin{lemma}\label{power>0}
For any positive real number $\alpha$, there exists a positive constant $c=c(\alpha)$, such that, for all $a,b \in \mathbb{R}$ 
\begin{equation}\label{ineqpositivepower}
    \frac{1}{c} \left| |b|^{\alpha-1} b - |a|^{\alpha-1} a \right| \leqslant \left( |a|+|b|\right)^{\alpha-1} |b-a| \leqslant c \left| |b|^{\alpha-1} b - |a|^{\alpha-1} a \right| \ .
\end{equation} 
\end{lemma}

\begin{lemma}\label{power>1}
For any real number $\gamma >1 $, there exists a positive constant $c=c(\gamma)$, such that, for all $a,b \in \mathbb{R}$ 
\begin{equation}\label{ineqpowerg1}
     |a-b|^{\gamma} \leqslant c \left| |a|^{\gamma-1} a - |b|^{\gamma-1} b\right| \ .
\end{equation} 
\end{lemma}

\noindent Next result refers to an anisotropic embedding and can be found in  \cite{DMV19} (for the elliptic anisotropic embedding we refer to \cite{Tro}).
\begin{lemma}\label{lem2.1}
Let $\alpha:=\frac{1}{N}\sum\limits_{i=1}^N \alpha_i$, $\alpha_i>0$, $\beta>0$, $p<N$ and $u \in  L^1(0, T; W^{1,1}_0(\Omega))$. Then there exists $\gamma >0$,
depending only on $N$, $\beta$, $p_1$, ..., $p_N$, $\alpha_1$, ..., $\alpha_N$ such that
\begin{equation*}
\iint\limits_{\Omega_T}|u|^{p_*}\,dx dt\leqslant \gamma \Big(\sup\limits_{0<t<T}\int\limits_{\Omega}|u(x, t)|^\beta dx\Big)^{\frac{p}{N}}
\prod\limits_{i=1}^N\Big(\iint\limits_{\Omega_T} |D_i u^{\alpha_i}|^{p_i}\,dx dt\Big)^{\frac{p}{N p_i}},\,\,p_*=\alpha p+\beta\frac{p}{N},
\end{equation*}
being the inequality trivial for unbounded right hand side.
\end{lemma}

\vspace{.2cm}

\noindent The following result corresponds to the energy estimates derived to \eqref{eq1.1} when the structure conditions \eqref{eq1.2} are satisfied (see \cite[Lemma 2.3]{CiaHenSavSkr2} for details). 
\begin{lemma}\label{lem2.3}
Let $u$ be a non-negative, local weak sub(super)-solution to  \eqref{eq1.1}-\eqref{eq1.2}. Then there exists a constant $\gamma >0$ depending only on the data $N$,  $p_1$, ..., $p_N$, $m_1$, ..., $m_N$, $K_{1}$, $K_{2}$, such that for every cylinder $Q_{\vec{r}, \eta}(y, \tau)\subset \Omega_T$, every $k\in \mathbb{R}_+$, and every piecewise smooth cutoff function $\zeta$, vanishing on $\partial K_{\vec{r}}(y)$, and such that $0\leqslant \zeta \leqslant 1$ there holds\[
\sup\limits_{\tau-\eta\leqslant t  \leqslant \tau}\int\limits_{K_{\vec{r}}(y)}g_{\pm}(u^{m_-}, k^{m_-})\,\zeta^{p_+}\,dx + \gamma^{-1}  \sum\limits_{i=1}^N\iint\limits_{Q_{\vec{r},\eta}(y, \tau)}u^{(m_i-m_-)(p_i-1)}|D_i (u^{m_-}-k^{m_-})_{\pm}|^{p_i}\zeta^{p_+}\,dx\,dt \]
\begin{eqnarray}\label{eq2.1}
&  \leqslant & \int\limits_{K_{\vec{r}}(y)\times\{\tau-\eta\}}g_{\pm}(u^{m_-}, k^{m_-})\,\,\zeta^{p_+}\,dx+\gamma\iint\limits_{Q_{\vec{r},\eta}(y, \tau)} g_{\pm}(u^{m_-}, k^{m_-})\,|\zeta_t|\,dx dt \nonumber \\ 
& & +\gamma \sum\limits_{i=1}^N\iint\limits_{Q_{\vec{r},\eta}(y, \tau)}u^{(m_i-m_-)(p_i-1)}(u^{m_-}-k^{m_-})_{\pm}^{p_i}|D_i \zeta|^{p_i}\,dx\,dt
\end{eqnarray}
\end{lemma} 
function $g_{\pm}$ is given by
\[g_{\pm}(u^{m}, k^{m}):=\int\limits^{(u^{m}-k^{m})_{\pm}}_0\big(k^{m}\pm s\big)^{\frac{1-m}{m}}\,s\,ds,\quad m>0,\]
and is bounded from below and above as follows (see Lemma 3.2 in \cite{BogDuzGiaLiaSch})
\begin{equation*}
\frac{1}{\gamma}\,(k^m+u^m)^{\frac{1}{m}-1}(u^m-k^m)^2_{\pm}\leqslant g_{\pm}(u^m, k^m)\leqslant \gamma\, (k^m+u^m)^{\frac{1}{m}-1}(u^m-k^m)^2_{\pm}, \quad \gamma=\gamma(m) \ .
\end{equation*}

\vspace{.2cm}

\noindent In what follows we present two De Giorgi type lemma, that can be found in full detail in \cite[Lemma 2.4]{CiaHenSavSkr2}. Here we are considering it written for $k=\xi\,\omega $, where 
$\mu^{\pm}$ and $\omega$ satisfy 
$$\mu^+\geqslant \sup\limits_{Q_{\vec{r}, \eta}(y, \tau)} u,\quad \mu^-\leqslant \inf\limits_{Q_{\vec{r}, \eta}(y, \tau)} u,\quad \omega\geqslant \mu^+-\mu^-.$$
\begin{lemma}\label{lem2.6}
Let $u$ be a non-negative, local weak super-solution to \eqref{eq1.1}, \eqref{eq1.2}. Fix $\xi$, $a\in (0, 1)$, then there exists $\nu_-\in (0, 1)$
depending only on the data, $a$, $\xi$, $\vec{r}$, $\eta$, $\omega$ such that if
\begin{equation}\label{eq2.5}
|Q_{\vec{r},\eta}(y, \tau)\cap\{u\leqslant \xi\,\omega\}|\leqslant \nu_- |Q_{\vec{r},\eta}(y, \tau)|,
\end{equation}
then
\begin{equation}\label{eq2.6}
u(x,t)\geqslant a \xi\,\omega,\quad a.e. \, (x,t)\in Q_{\frac{1}{2}\vec{r}, \frac{1}{2}\eta}(y, \tau).
\end{equation}
\end{lemma}

\begin{remark}
In what follows we will use Lemma \ref{lem2.6} with the choices
$$\eta= b\,r^{p_+}\,\omega^{1-\lambda_+},\quad r_i=r^{\frac{p_+}{p_i}}\omega^{-\frac{\lambda_+-\lambda_i}{p_i}},\,\,i=1, ...,N,$$
with some $b\in (0, 1)$, $r>0$. Simple calculations give that the number
$\displaystyle{\nu:=\frac{1}{\gamma}\,b^{\frac{N}{p}}}$ satisfies $\nu\leqslant \nu_-.$ 
\end{remark}

\vspace{.2cm}

\noindent The final results concern the expansion of positivity for the local weak super-solutions to \eqref{eq1.1}-\eqref{eq1.2} in $\Omega_T$ in the singular case $\lambda_+<1$ and in the degenerate case $\lambda_+>1$, respectively (for more details we refer to \cite{CiaHenSavSkr2}). 

\begin{lemma}\label{lem3.2}
Let $u$ be a non-negative, local weak super-solution to \eqref{eq1.1}-\eqref{eq1.2} in $\Omega_T$ and $\lambda_+<1$. If
\begin{equation}\label{eq3.2}
|K^{k}_r(\bar{x})\cap\{u(\cdot, \bar{t})\geqslant k\}|\geqslant \alpha_0 |K^{k}_r(\bar{x})|,
\end{equation}
with some positive $k$, $r>0$ and $\alpha_0 \in (0, 1)$, then there exist constants $\delta$, $\varepsilon$, $\epsilon_{*}\in(0,1)$ and $C_{*}>1$, depending only on the data and $\alpha_0$, such that
\begin{equation}\label{eq3.3}
u(x, t)\geqslant \frac{k}{C_{*}},\quad a.e. \ \  x\in K^{k/C_{*}}_{r}(\bar{x}),
\end{equation}
and for all $t$
\begin{equation}\label{eq3.4}
\bar{t}+(1-\varepsilon) \delta \,r^{p_{+}}\,\,k^{1-\lambda_{+}} \leqslant t\leqslant \bar{t}+ \delta \,r^{p_{+}}\,\,k^{1-\lambda_{+}},
\end{equation}
provided that
\begin{equation}\label{eq3.5}
\lambda_+-\lambda_-\leqslant \epsilon_*.
\end{equation}
\end{lemma}

\begin{lemma}\label{lem4.1}
Let $u$ be a non-negative, local weak super-solution to \eqref{eq1.1}-\eqref{eq1.2} in $\Omega_T$ and $\lambda_+>1$. If
\begin{equation}\label{eq4.1}
|K^{k}_r(\bar{x})\cap\{u(\cdot, \bar{t})\geqslant k\}|\geqslant \alpha_0 |K^{k}_r(\bar{x})|,
\end{equation}
with some positive $k$, $r>0$ and $\alpha_0 \in (0, 1)$, then there exist constants  $\epsilon_{*}\in(0,1)$, $B>1$ and $C_{*}>1$, depending only on the data and $\alpha_0$, such that
\begin{equation}\label{eq4.2}
u(x, t)\geqslant \frac{k}{C_{*}},\quad a.e. \ \  x\in K^{k/C_{*}}_{r}(\bar{x}),
\end{equation}
and for all $t$
\begin{equation}\label{eq4.3}
\bar{t}+ \frac{1}{2} B \,r^{p_{+}}\,\,k^{1-\lambda_{+}} \leqslant t\leqslant \bar{t}+ B \,r^{p_{+}}\,\,k^{1-\lambda_{+}},
\end{equation}
provided that
\begin{equation}\label{eq4.4}
\lambda_+-\lambda_-\leqslant \epsilon_*.
\end{equation}
\end{lemma}

\subsection{De Giorgi type lemmas}\label{SS:DEGlemmas}

\noindent The lemmas in this  section are well known in the standard (isotropic) case as De Giorgi type lemmas, here we show that similar results are valid in the anisotropic case.  

\noindent Construct the cylinder $Q_{2\vec{r}, 2\eta}(y, \tau)\subset \Omega_T$, denote by $\mu^{\pm}$ and $\omega$ the following numbers
$$\mu^+\geqslant \sup\limits_{Q_{\vec{r}, \eta}(y, \tau)} u,\quad \mu^-\leqslant \inf\limits_{Q_{\vec{r}, \eta}(y, \tau)} u,\quad \omega\geqslant \mu^+-\mu^-.$$

\begin{lemma}\label{lem2.7}
Let $u$ be a non-negative, locally bounded, local weak sub-solution to \eqref{eq1.1}-\eqref{eq1.2}. Then there exists $\nu_+ \in (0, 1)$ depending only on the data, $\vec{r}$, $\eta$, $\omega$ such that if
\begin{equation}\label{eq2.8}
|Q_{\vec{r},\eta}(y, \tau)\cap\{u\geqslant \mu^{+}-\frac{1}{4}\omega\big\}|\leqslant \nu_+ |Q_{\vec{r},\eta}(y, \tau)|,
\end{equation}
then
\begin{equation}\label{eq2.9}
u(x, t)\leqslant \mu^+- \frac{1}{8} \omega,\quad a.e. \ \  (x,t)\in Q_{\frac{1}{2}\vec{r}, \frac{1}{2}\eta}(y, \tau),
\end{equation}
provided that
\begin{equation}\label{eq2.10}
\frac{1}{2}\,\omega\leqslant \mu^+\leqslant 2\,\omega.
\end{equation}
\end{lemma}

\begin{proof} For $j=0,1,\cdots$, consider 
\[ k_j=\mu^{+}-\frac{1}{8}\omega(1+2^{-j}), \qquad   \vec{r}_j:=\frac{1}{2}\vec{r}(1+2^{-j}), \qquad  
\eta_j:=\frac{1}{2}\eta (1+2^{-j}) \ , \]  the cubes $K_j:=K_{\vec{r}_j}(y)$ and cylinders  $Q_j:=K_j\times I_j$, where $I_j:= (\tau-\eta_j, \tau)$ and cutoff functions $\zeta^{(1)}_j(x) \times \zeta^{(2)}_j(t)$, defined in $Q_j$ and satisfying 
 $\zeta^{(1)}_j(x)\in C^1_0(K_j)$, $\zeta^{(1)}_j(x)=1$ in $K_{j+1}$, $0\leqslant \zeta^{(1)}_j(x)\leqslant 1$ and $|D_i \zeta^{(1)}_j(x)|\leqslant \frac{\gamma 2^{j}}{r_i}$, $i=1, ..., N$ and $\zeta^{(2)}_j(t)\in C^1(\mathbb{R})$, $\zeta^{(2)}_j(t)=0$ for $t\leqslant \tau-\eta_j$, $\zeta^{(2)}_j(t)=1$ for 
$t\geqslant \tau-\eta_{j+1}$, $0\leqslant \zeta^{(2)}_j(t)\leqslant 1$, $|\zeta^{(2)}_{j, t}|\leqslant \frac{\gamma 2^j}{\eta}$. Under this setting, write \eqref{eq2.1} for $\varphi(x,t)= (u^{m_{-}}-k^{m_{-}}_j)_{+} \zeta_j^{p_+} $, defined in $Q_j$, 
\begin{equation}\label{est}
\begin{aligned}
\sup\limits_{t \in I_j}&\int\limits_{K_j \times\{t\}}g_{+}(u^{m_{-}}, k^{m_{-}}_j)\,\zeta^{p_{+}}_j\,dx+\gamma^{-1}\sum\limits_{i=1}^N \iint\limits_{Q_j}u^{(m_i-m_{-})(p_i-1)}\left| D_i(u^{m_{-}}-k^{m_{-}}_j)_{+} \right|^{p_{i}}\,\zeta^{p_{+}}_j dx\,dt \\
&  \leqslant  \frac{\gamma 2^{j}}{\eta}\iint\limits_{Q_j} g_{+}(u^{m_{-}}, k_j^{m_{-}})\,\,dx dt +\gamma \sum\limits_{i=1}^{N} \frac{2^{j\ p_i}}{r^{p_i}_i} \iint\limits_{Q_j}u^{(m_i-m_{-})(p_i-1)}(u^{m_{-}}-k^{m_{-}}_j)_{+}^{p_{i}}\,dx\,dt  \ .
\end{aligned}
\end{equation}
Recalling that 
\[\frac{1}{\gamma(m_-)} \omega^{1-m_-} \;  (u^{m_{-}}- k_j^{m_{-}})_+^2 \leqslant g_{+}(u^{m_{-}}, k_j^{m_{-}})\leqslant \gamma(m_-) \; \omega^{1-m_-}  (u^{m_{-}}- k_j^{m_{-}})_+^2  \]
\[  k_j <u \leqslant 2 \omega , \qquad  \dfrac{1}{4} \omega  \leqslant k_j \leqslant \dfrac{15}{8}\omega\]
and considering $\alpha_i= \frac{\lambda_i+m_{-}}{m_- p_i}>0$, from the previous estimate we arrive at 
\begin{multline}
\omega^{1-m_{-}}\sup\limits_{t\in I_j}\int\limits_{K_j\times\{t\}}\negthickspace\negthickspace(u^{m_{-}}-k^{m_{-}}_j)_{+}^2\,\zeta^{p_{+}}_j\,dx+ \gamma^{-1} \sum\limits_{i=1}^N \iint\limits_{Q_j}\left|D_i(u^{m_{-}}-k^{m_{-}}_j)^{\alpha_i}_{+} \right|^{p_{i}} \negthickspace\negthickspace \zeta^{p_{+}}_j\negthickspace dx\,dt\\
\leqslant \gamma 2^{j \gamma}\ \omega^{\lambda_{+}+m^{-}}\Big(\frac{\omega^{1-\lambda_{+}}}{\eta}+\sum\limits_{i=1}^N \frac{\omega^{\lambda_i-\lambda_{+}}}{r^{p_{i}}_i}\Big)|A_j|, 
\end{multline}
where $A_j:=Q_j\cap [u > k_j]$.

\noindent On the one hand, by applying H\"{o}lder's inequality together with the Sobolev embedding Lemma~\ref{lem2.1}, with $\alpha_i=\frac{\lambda_i +m^{-}}{m^{-}\,p_i}$, $\alpha=\frac{1}{N}\sum\limits_{i=1}^n\alpha_i$, and $\beta=2$, we get
\begin{eqnarray*}
\iint\limits_{Q_j}\left((u^{m_{-}}-k^{m_{-}}_j)_{+}\zeta^{p_{+}}_j\right)^{\alpha p}\,dxdt & \leqslant & \Big(\iint\limits_{Q_j}\big[(u^{m_{-}}-k^{m_{-}}_j)_{+}\zeta^{p_{+}}_j\big]^{p\frac{\alpha N+2}{N}}\,dxdt\Big)^{\frac{N}{N+2/\alpha}}|A_j|^{\frac{2/\alpha}{N+2/\alpha}}\\ & \leqslant &
\gamma\Big(\sup\limits_{t\in I_j}\int\limits_{K_j\times\{t\}}(u^{m_{-}}-k_j^{m_{-}})^{2}_{+}\zeta^{p_{+}}_j\,dx\Big)^{\frac{p}{N+2/\alpha}} \\
& & \times \Big(\sum\limits_{i=1}^N \iint\limits_{Q_j}\left|D_i\left((u^{m_{-}}-k^{m_{-}}_j)^{\alpha_i}_{+}\zeta^{p_{+}}_j\right) \right|^{p_{i}}\,dxdt\Big)^{\frac{N}{N+2/\alpha}}|A_j|^{\frac{2/\alpha}{N+2/\alpha}} \\
 & \leqslant & \gamma 2^{j\gamma} \omega^{\frac{(m_{-}-1)p}{N+2/\alpha}+\frac{(\lambda_{+}+m_{-})(N+p)}{N+2/\alpha}}\\
    & & \times\Big(\frac{\omega^{1-\lambda_{+}}}{\eta}+\sum\limits_{i=1}^N \frac{\omega^{\lambda_i -\lambda_{+}}}{r^{p_{i}}_i}\Big)^{\frac{N+p}{N+2/\alpha}}\times |A_j|^{1+\frac{p}{N+2/\alpha}}.
\end{eqnarray*}
On the other hand, by recalling \eqref{eq2.10} and \eqref{ineqpositivepower},
\[
\iint\limits_{Q_j}(u^{m_{-}}-k^{m_{-}}_j)_{+}^{\alpha p}\zeta^{p_{+}\alpha p}_j\,dxdt \geqslant 1(k^{m_{-}}_{j+1}-k^{m_{-}}_j)^{\alpha p}
|A_{j+1}| \geqslant \gamma \, 
\omega^{m^{-} \alpha p} \, 2^{- j\alpha p} |A_{j+1}|. 
\]
Gathering the previous estimates and setting $y_j:=\dfrac{|A_j|}{|Q_j|}$, we get
\[
y_{j+1}\leqslant 
\gamma 2^{j\gamma} \Big(\frac{\eta}{\omega^{1-\lambda_{+}}}\Big)^{\frac{p}{N+2/\alpha}}\Big(\prod\limits_{i=1}^N \frac{r_i}{\omega^{\frac{\lambda_i-\lambda_{+}}{p_i}}}\Big)^{\frac{p}{N+2/\alpha}}\Big(\frac{\omega^{1-\lambda_{+}}}{\eta}+\sum\limits_{i=1}^N \frac{\omega^{\lambda_i-\lambda_{+}}}{r^{p_{i}}_i}\Big)^{\frac{N+p}{N+2/\alpha}}y_{j}^{1+\frac{p}{N+2/\alpha}}.\]
The proof now follows once we take $\nu_+$ such that
\begin{equation}\label{eq2.8}
\nu_+:=\frac{1}{\gamma}\frac{\omega^{1-\lambda_{+}}}{\eta}\Bigg(\prod\limits_{i=1}^N \frac{r_i}{\omega^{\frac{\lambda_i-\lambda_{+}}{p_i}}}\Bigg)^{-1}\Big(\frac{\omega^{1-\lambda_{+}}}{\eta}+\sum\limits_{i=1}^N \frac{\omega^{\lambda_i-\lambda_{+}}}{r^{p_{i}}_i}\Big)^{-\frac{N+p}{p}}.
\end{equation}
\end{proof}

\begin{remark}
In what follows we will use Lemma \ref{lem2.7} with the choices
$$\eta= b\,r^{p_+}\,\omega^{1-\lambda_+},\quad r_i=r^{\frac{p_+}{p_i}}\omega^{-\frac{\lambda_+-\lambda_i}{p_i}},\,\,i=1, ...,N,$$
with some $b\in (0, 1]$, $r>0$. One can easy show that 
$\displaystyle{\nu:=\frac{1}{\gamma}\,b^{\frac{N}{p}}}$ satisfies $\nu\leqslant \nu_+.$ 
\end{remark}

\vspace{.2cm}
\noindent The two following results are presented and proved within the degenerate range $\lambda_+>1$, distinguishing between $p_+>2$ and $p_+< 2$. They state that if a non-negative, locally bounded, local weak sub-solution to \eqref{eq1.1} satisfy an initial data then something alike is also valid in a proper cylinder; this kind of results is known in the literature as De Giorgi type lemma involving initial data.

\begin{lemma}\label{lem2.8}
Let $u$ be a non-negative, locally bounded, local weak sub-solution to \eqref{eq1.1}-\eqref{eq1.2} in $\Omega_T$, $\lambda_+>1$ and $p_+>2$.
Assume that, for some $\xi_0\in (0, \frac{1}{2})$,
\begin{equation}\label{eq2.13}
u(x, \tau-\eta)\leqslant \mu^+- \xi_0\,\omega,\quad a.e. \, x\in K_{\vec{r}}(y).
\end{equation}
Then there exists  $\xi\in (0, \xi_0)$, depending only on the data,    $\vec{r}$, $\eta$ and $\omega$, such that 
\begin{equation}\label{eq2.14}
u(x, t)\leqslant \mu^+- \frac{1}{4}\xi \omega,\quad  a.e. \ (x, t)\in Q_{\frac{1}{2}\vec{r}, \eta}(y, \tau),
\end{equation}
provided that
\begin{equation}\label{eq2.15}
\frac{1}{2}\,\omega \leqslant \mu^+\leqslant 2\,\omega,
\end{equation}
and
\begin{equation}\label{eq2.16}
 p_+-p_-\leqslant \epsilon^{(1)}_*
\end{equation}
for some sufficiently small $\epsilon^{(1)}_*$ depending only on the data and $\xi$.
\end{lemma}

\begin{proof}
The proof is divided into two parts: first and under the assumption that \eqref{eq2.13} holds, we prove \\
\textbf{Claim.} For any $\nu \in (0, 1)$ there exists $\xi\in (0, \xi_0)$ depending only on the data, $\eta$, $\vec{r}$, $\omega$ and $\nu$
such that
\begin{equation}\label{eq2.17}
|K_{\frac{3}{4}\vec{r}}(y)\cap\{u(\cdot, t)\geqslant \mu^+-\frac{1}{2}\xi\,\omega\}|\leqslant \nu |K_{\frac{3}{4}\vec{r}}(y)|,
\end{equation}
for all $\tau-\eta\leqslant t\leqslant \tau$.

\vspace{.2cm}

\noindent In fact, consider the estimates \eqref{eq2.1} written over the cylinder $Q_{\vec{r}, \eta}(y, \tau)$ for $k=\mu^+-\xi\,\omega$, and a time independent cutoff function $\zeta(x)\in C^1_0(K_{\vec{r}}(y))$ satisfying $\zeta(x)=1$ in $K_{\frac{3}{4}\vec{r}}(y)$ and  $0\leqslant \zeta(x)\leqslant 1$ and $|D_i \zeta(x)|\leqslant \frac{\gamma}{r_i}$, $i=1, ..., N$. Then, for all $\tau-\eta \leqslant t \leqslant \tau$, on the one hand we get 
\begin{eqnarray*}
    \omega^{1-m_-} \int\limits_{K_{\vec{r}}(y)}(u^{m_-}-k^{m_-})^2_+\,\zeta^{p_+}\,dx & \leqslant &
\gamma \sum\limits_{i=1}^N\frac{1}{r^{p_i}_i}\omega^{(m_i-m_-)(p_i-1)}\iint\limits_{Q_{\vec{r}, \eta}(y, \tau)}(u^{m_-}-k^{m_-})^{p_i}_+\,dx dt\\
& \leqslant & \gamma \eta\,|K_{\vec{r}}(y)| \,\sum\limits_{i=1}^N \frac{1}{r^{p_i}_i}\omega^{\lambda_i+m_-}\,\xi^{ p_i}
\end{eqnarray*}
due to estimates \eqref{eq2.13} and \eqref{eq2.15} and observing that, by using \eqref{ineqpositivepower},
\[ (u^{m_-}-k^{m_-})_+ \leqslant {\mu^+}^{m_-} -k^{m_-} \leqslant \gamma(m_-) ( \mu^++k)^{m_--1} (\mu^+-k) \leqslant \gamma(m_-) \xi \, \omega^{m_-} \ . \]
On the other hand
\[\int\limits_{K_{\vec{r}}(y)}(u^{m_-}-k^{m_-})^2_+\,\zeta^{p_+}\,dx \geqslant \frac{1}{\gamma} \left( \xi \, \omega^{m_-}\right)^2\,|K_{\frac{3}{4}\vec{r}}(y)\cap\{u(\cdot, t) > \mu^+-\frac{1}{2}\xi\,\omega\}| \ . \]
The required \eqref{eq2.17} follows, provided that
$\xi$ and $\epsilon^{(1)}_*$ are chosen to satisfy
\begin{equation}\label{eq2.18}
\gamma\,\xi^{p_+-2}\,\eta\,\sum\limits_{i=1}^N\frac{\omega^{\lambda_i-1}}{r^{p_i}_i}\,\leqslant \nu,
\end{equation}
and 
\[\xi^{-(p_+-p_-)}\leqslant \xi^{-\epsilon^{(1)}_*}\leqslant 2 \ .\]

\noindent The next and final step in the proof consists in fixing $\nu$ (thereby choosing $\xi$) verifying \eqref{eq2.18} so that a fast geometric convergence lemma can be applied and \eqref{eq2.14} be derived. For that purpose, and for $j=0,1,\cdots$, choose $k_j:=\mu^+- \frac{1}{4}\xi \omega (1+2^{-j})$,  $\vec{r}_j:=\frac{\vec{r}}{2}(1+2^{-j-1})$, construct the cubes
$K_j:=K_{\vec{r}_j}(y)$ and the cylinders $Q_j:=K_j\times I$, $I:=(\tau -\eta, \tau)$, and define cutoff functions $\zeta_j(x)\in C^1_0(K_j)$, $\zeta_j(x)=1$ in $K_{j+1}$, $0\leqslant \zeta_j(x)\leqslant 1$,  $|D_i \zeta_j(x)|\leqslant \frac{\gamma 2^j}{r_i}$, $i=1, ..., N$. Under this setting, we apply estimate \eqref{eq2.1} and then obtain
\[ \omega^{1-m_-} \sup\limits_{t\in I}\int\limits_{K_j\times\{t\}}(u^{m_-}-k_j^{m_-})^2_{+}\zeta^{p_{+}}_j\,dx+ \gamma^{-1} \sum\limits_{i=1}^{N} \omega^{(m_i-m_-)(p_i-1)}\iint\limits_{Q_j}|D_i(u^{m_-}-k_j^{m_-})_{+}|^{p_{i}}\zeta^{p_{+}}_j\,dx\,dt \]
\vspace{-.5cm}
\begin{eqnarray*}
&\leqslant &  \gamma 2^{j\gamma}\sum\limits_{i=1}^N \frac{\omega^{\lambda_i+m_-}}{r^{p_i}_i}\xi^{p_i}|Q_j\cap\{u > k_j\}|\leqslant \gamma 2^{j\gamma}\xi^{p_-} \omega^{\lambda_++m_-} \sum\limits_{i=1}^N \frac{\omega^{\lambda_i-\lambda_+}}{r^{p_i}_i}|Q_j\cap\{u > k_j\}| \ .
\end{eqnarray*}
Proceeding analogously to Lemma \ref{lem2.7} (but now for $\alpha_i=1$, for all $i=1, \cdots,N$, and therefore $\alpha=1$), we get
\begin{eqnarray*}
    |Q_{j+1}\cap\{u > k_{j+1}\}| & \leqslant  &
\gamma 2^{j\gamma}\,\omega^{\frac{p}{N+2}(\sum_i\frac{\lambda_+-\lambda_i}{p_i} + \lambda_+-1)} \, \xi^{p_-\frac{p+N}{N+2} - p} \, \left(\sum\limits_{i=1}^N \frac{ \omega^{\lambda_i-\lambda_+}}{r^{p_i}_i}\right)^{\frac{N+p}{N+2}}\\
& & \times |Q_j\cap\{u > k_j\}|^{1+\frac{p}{N+2}}.
\end{eqnarray*} 
and then, for $y_j:=\dfrac{|Q_j\cap\{u > k_j\}|}{|Q_j|}$, we obtain the algebraic recursive relation
\begin{eqnarray}\label{eq2.19}
y_{j+1} &\leqslant & \gamma 2^{j\gamma}\,\xi^{\frac{p_-(N+p)}{N+2}- p}\,\Big(\frac{\eta}{\omega^{1-\lambda_{+}}}\Big)^{\frac{p}{N+2}}\bigg(\prod_{i=1}^N \frac{r_i}{ \omega^{\frac{\lambda_i-\lambda_{+}}{p_i}}}\bigg)^{\frac{p}{N+2}}\Big(\sum\limits_{i=1}^N \frac{ \omega^{\lambda_i-\lambda_+}}{r^{p_i}_i}\Big)^{\frac{N+p}{N+2}} y_j^{1+\frac{p}{N+2}} \nonumber \\ 
& \leqslant & \gamma 2^{j\gamma}\,\Big(\frac{\eta}{\omega^{1-\lambda_{+}}}\Big)^{\frac{p}{N+2}}\bigg(\prod_{i=1}^N \frac{r_i}{ \omega^{\frac{\lambda_i-\lambda_{+}}{p_i}}}\bigg)^{\frac{p}{N+2}}\Big(\sum\limits_{i=1}^N \frac{ \omega^{\lambda_i-\lambda_+}}{r^{p_i}_i}\Big)^{\frac{N+p}{N+2}}
y_j^{1+\frac{p}{N+2}}
\end{eqnarray}
since $p_-(N+p)-p(N+2)=(p_--p)N+p(p_--2)\geqslant -(p_+-p_-)(N+p)\geqslant -\epsilon^{(1)}_*(N+p)$ and  $\epsilon^{(1)}_*$ chosen as before. 

\noindent The previous recursive algebraic inequality gives us \eqref{eq2.14}, provided that
\[ y_0\leqslant \frac{1}{\gamma}\frac{\omega^{1-\lambda_+}}{\eta}\bigg(\prod_{i=1}^N \frac{r_i}{ \omega^{\frac{\lambda_i-\lambda_{+}}{p_i}}}\bigg)^{-1}\Big(\sum\limits_{i=1}^N \frac{ \omega^{\lambda_i-\lambda_+}}{r^{p_i}_i}\Big)^{-\frac{N+p}{p}}. \]
By \eqref{eq2.18}, this inequality holds if $\xi$ is chosen to satisfy
\begin{equation}\label{eq2.20}
\xi^{p_+-2}\leqslant \frac{1}{\gamma}\Big(\frac{\omega^{1-\lambda_+}}{\eta}\Big)^2\bigg(\prod_{i=1}^N \frac{r_i}{ \omega^{\frac{\lambda_i-\lambda_{+}}{p_i}}}\bigg)^{-1}
\Big(\sum\limits_{i=1}^N \frac{ \omega^{\lambda_i-\lambda_+}}{r^{p_i}_i}\Big)^{-\frac{N+2 p}{p}}.
\end{equation}

\end{proof}

\begin{remark}
In what follows we will use Lemma \ref{lem2.8} with the choices
$$\eta= B\,r^{p_+}\,\omega^{1-\lambda_+},\quad r_i=r^{\frac{p_+}{p_i}}\,\omega^{-\frac{\lambda_+-\lambda_i}{p_i}},\,\,i=1, ...,N,$$
with some $r>0$ and $B\gg 1.$ One can easily prove that 
$\displaystyle{\xi =\Big(\frac{1}{\gamma B}\Big)^{\frac{2}{p_+-2}}}$
satisfies \eqref{eq2.20}. In this case the  $\epsilon^{(1)}_*$, from \eqref{eq2.16}, depends only on the data and $B$.
\end{remark}

\noindent The next result is a variant of De Giorgi type lemma involving initial data and it is presented within the range $\lambda_+>1$ and $p_+<2$.

\begin{lemma}\label{lem2.9}
Let $u$ be a non negative, locally bounded, local weak sub-solution to \eqref{eq1.1}-\eqref{eq1.2} in $\Omega_T$, $\lambda_+>1$ and $p_+ <2$. Assume that, for some $\xi_0\in (0, 1/4]$, there holds
\begin{equation}\label{eq2.21}
u(x, \tau)\leqslant \mu^+- \xi_0\,\omega,\quad a.e.\  \ x\in K_{\vec{r}}(y).
\end{equation}
Then there exists $\bar{\nu}_0 \in (0, 1)$, depending only on the data, $\vec{r}$, $\omega$  and $\xi_0$, such that 
\begin{equation}\label{eq2.22}
u(x, t)\leqslant \mu^+ - \frac{1}{4}\xi_0 \omega,\quad  a.e. \ (x, t)\in K_{\frac{1}{4}\vec{r}}(y)\times(\tau, \tau +\bar{\nu}_0),
\end{equation}
provided that
\begin{equation}\label{eq2.23}
\frac{1}{2}\,\omega\leqslant \mu^+\leqslant 2\,\omega \ ,
\end{equation}
and 
\begin{equation}\label{eq2.16}
 p_+-p_-\leqslant \epsilon^{(1)}_*
\end{equation}
for some sufficiently small $\epsilon^{(1)}_*$ depending only on the data and $\xi_0$.
\end{lemma}

\begin{proof}
The structure of the proof is quite similar to the one presented previously for $p_+>2$: assuming \eqref{eq2.21} holds, we start to prove a claim from which we can obtain \eqref{eq2.22}.

\textbf{Claim.} For any $\nu \in (0, 1)$ there exists $\bar{\nu}_0$, depending only on the data, $\nu$, $\vec{r}$, $\omega$ and $\xi_0$, such that
\begin{equation}\label{eq2.25}
|K_{\frac{1}{2}\vec{r}}(y)\cap\{u(\cdot, t)\geqslant \mu^+-\frac{1}{2}\,\xi_0\,\omega\}|\leqslant \nu |K_{\frac{1}{2}\vec{r}}(y)|, \quad \forall t\in (\tau, \tau +\bar{\nu}_0) \ .
\end{equation}
\noindent Indeed, by considering time-independent cutoff function  $\zeta(x)\in C^1_0(K_{\vec{r}}(y))$, $0\leqslant \zeta(x)\leqslant 1$, $\zeta(x)=1$ in $K_{\frac{1}{2}\vec{r}}(y)$, and $|D_i \zeta(x)|\leqslant \frac{\gamma}{r_i}$, $i=1, ..., N$, recalling that \eqref{eq2.21} and \eqref{eq2.23} hold, observing that 
\[ \frac{1}{4} \omega \leqslant (\frac{1}{2}- \xi_0)\omega <  u\leqslant 2\omega, \quad \text{on the set}\quad u > k=\mu^+-\xi_0\,\omega \ , \]
from \eqref{eq2.1} we obtain, for all $\tau\leqslant t \leqslant \tau+\bar{\nu}_0$,
\begin{eqnarray*}
    \int\limits_{K_{\vec{r}}(y)\times\{t\}}(u^{m_-}-k^{m_-})^2_+\,\zeta^{p_+}\,dx & \leqslant &
 \gamma \bar{\nu}_0\,|K_{\vec{r}}(y)| \,\sum\limits_{i=1}^N \frac{1}{r^{p_i}_i}\ \xi_0^{p_i} \, \omega^{\lambda_i-2m_--1}  \\
 & = &
 \gamma \, \bar{\nu}_0\, \xi_0^{p_+}\, |K_{\vec{r}}(y)| \, \sum\limits_{i=1}^N \xi_0^{p_i-p_+}\frac{\omega^{\lambda_i-2m_--1}}{r^{p_i}_i}  \ ,
\end{eqnarray*} 
 on the other hand and for all $\tau\leqslant t \leqslant \tau+\bar{\nu}_0$,
\begin{eqnarray*} 
\int\limits_{K_{\vec{r}}(y) \times \{t\}}(u^{m_-}-k^{m_-})^2_+\,\zeta^{p_+}\,dx & \geqslant & \int\limits_{K_{\frac{1}{2}\vec{r}}(y) \cap  \{u(\cdot, t) > \mu^+-\frac{1}{2}\xi_0\,\omega\} }(u^{m_-}-k^{m_-})^2_+ \,dx  \\
& \geq & \left(\frac{1}{2} \xi_0\right)^2 \omega^{2m_-} \, |K_{\frac{1}{2}\vec{r}}(y)\cap\{u(\cdot, t) > \mu^+-\frac{1}{2}\xi_0\,\omega\}|  \ . 
 \end{eqnarray*}
Hence, for all $\tau\leqslant t \leqslant \tau+\bar{\nu}_0$,
\[ |K_{\frac{1}{2}\vec{r}}(y)\cap\{u(\cdot, t) > \mu^+-\frac{1}{2}\xi_0\,\omega\}|\leqslant  \gamma \, \bar{\nu}_0 \, \xi^{p_+-2}_0 \,\sum\limits_{i=1}^N \xi_0^{p_i-p_+}\, \frac{\omega^{\lambda_i-1}}{r^{p_i}_i}\, |K_{\vec{r}}(y)| \ .
\]
The required \eqref{eq2.25} follows, provided that
$\bar{\nu}_0$ satisfies
\begin{equation}\label{eq2.26}
\gamma \, \bar{\nu}_0 \, \xi^{p_+ -2}_0 \,\sum\limits_{i=1}^N \frac{\omega^{\lambda_i-1}}{r^{p_i}_i} \leqslant \nu \ ,
\end{equation}
and $\epsilon_*^{(1)}$ is such that 
\[ \xi_0^{p_--p_+} \leqslant \xi_0^{-\epsilon_*^{(1)}} \leqslant 2 \ . \]
Being the claim proved, we now proceed with the remainder of the proof. For $j=0,1, \cdots$, consider the sequence of levels $ k_j:=\mu^+-\frac{1}{4}\xi_0 \omega-
\frac{1}{4}\xi_0 \omega 2^{-j}$, define  
\[ \vec{r}_j:=\frac{\vec{r}}{4}(1+2^{-j}) \ , \quad 
K_j:=K_{\vec{r}_j}(y) \ , \quad Q_j:=K_j\times I , \quad I:=(\tau, \tau +\bar{\nu}_0) \ \] and consider cutoff functions $\zeta_j(x)\in C^1_0(K_j)$, $0\leqslant \zeta_j(x)\leqslant 1$, $\zeta_j(x)=1$ in $K_{j+1}$,  $|D_i \zeta_j(x)|\leqslant \frac{\gamma 2^j}{r_i}$, $i=1, \cdots, N$. Under this setting and recalling \eqref{eq2.21} and \eqref{eq2.23},  estimate \eqref{eq2.1} reads
\[ \omega^{1-m_-} \sup\limits_{t\in I}\int\limits_{K_j\times\{t\}}(u^{m_-}-k_j^{m_-})^2_{+}\zeta^{p_{+}}_j\,dx+ \gamma^{-1} \sum\limits_{i=1}^{N} \omega^{(m_i-m_-)(p_i-1)}\iint\limits_{Q_j}|D_i(u^{m_-}-k_j^{m_-})_{+}|^{p_{i}}\zeta^{p_{+}}_j\,dx\,dt \]
\vspace{-.5cm}
\begin{eqnarray*}
&\leqslant &  \gamma 2^{j\gamma}\sum\limits_{i=1}^N \frac{\omega^{\lambda_i+m_-}}{r^{p_i}_i}\xi_0^{p_i}|Q_j\cap\{u > k_j\}|\leqslant \gamma 2^{j\gamma}\xi_0^{p_+} \omega^{\lambda_++m_-} \sum\limits_{i=1}^N \frac{\omega^{\lambda_i-\lambda_+}}{r^{p_i}_i}|Q_j\cap\{u > k_j\}| \ ,
\end{eqnarray*}
recalling the choice of $\epsilon_*^{(1)}$. From here the proof brings no extra difficulties and we proceed analogously to the proof of Lemma 
\ref{lem2.7}, therefore obtaining, for $y_j:=\frac{|Q_j\cap\{u\geqslant k_j\}|}{|Q_j|}$, 
\begin{eqnarray*}
y_{j+1} &\leqslant &\gamma 2^{j\gamma}\,\xi^{-p+ p_+\frac{N+p}{N+2}}_0\,\Big(\frac{\bar{\nu}_0}{\omega^{1-\lambda_{+}}}\Big)^{\frac{p}{N+2}}\bigg(\prod_{i=1}^N \frac{r_i}{ \omega^{\frac{\lambda_i-\lambda_{+}}{p_i}}}\bigg)^{\frac{p}{N+2}} \Big(\sum\limits_{i=1}^N \frac{ \omega^{\lambda_i-\lambda_+}}{r^{p_i}_i}\Big)^{\frac{N+p}{N+2}}
y_j^{1+\frac{p}{N+2}}\\
&\leqslant &\gamma 2^{j\gamma}\,\xi^{(p_+-2)\frac{p}{N+2}}_0\,\Big(\frac{\bar{\nu}_0}{\omega^{1-\lambda_{+}}}\Big)^{\frac{p}{N+2}}\bigg(\prod_{i=1}^N \frac{r_i}{ \omega^{\frac{\lambda_i-\lambda_{+}}{p_i}}}\bigg)^{\frac{p}{N+2}} \Big(\sum\limits_{i=1}^N \frac{ \omega^{\lambda_i-\lambda_+}}{r^{p_i}_i}\Big)^{\frac{N+p}{N+2}}
y_j^{1+\frac{p}{N+2}}
\end{eqnarray*}
which yields $\lim\limits_{j\rightarrow \infty} y_j=0$, provided that
\begin{equation*}
y_0\leqslant \frac{1}{\gamma}\xi^{2-p_{+}}_0\,\frac{\omega^{1-\lambda_+}}{\bar{\nu}_0}\bigg(\prod_{i=1}^N \frac{r_i}{ \omega^{\frac{\lambda_i-\lambda_{+}}{p_i}}}\bigg)^{-1} \Big(\sum\limits_{i=1}^N \frac{ \omega^{\lambda_i-\lambda_+}}{r^{p_i}_i}\Big)^{-\frac{N+p}{p}}.
\end{equation*}
By \eqref{eq2.26} this inequality holds, if $\bar{\nu}_0$ is chosen to satisfy
\begin{equation}\label{eq2.27}
\bar{\nu}_0\leqslant \frac{1}{\gamma}\xi^{2-p_+}_0\,\omega^{1-\lambda_+}\bigg(\prod_{i=1}^N \frac{r_i}{ \omega^{\frac{\lambda_i-\lambda_{+}}{p_i}}}\bigg)^{-\frac{1}{2}} \Big(\sum\limits_{i=1}^N \frac{ \omega^{\lambda_i-\lambda_+}}{r^{p_i}_i}\Big)^{-\frac{N+2p}{2p}}.
\end{equation}
\end{proof}


\subsection{Integral estimates leading to Harnack-type inequalities}

\begin{lemma}\label{lem2.4}
Let $u$ be a non-negative, local weak super-solution to \eqref{eq1.1}-\eqref{eq1.2} in $\Omega_T$ and fix $\alpha>0 $ and $1<\beta <1+\frac{1}{m_-}$. Then there exists $\gamma >0$ depending only on the data and $\beta$ such that for every $\varepsilon>0$, every cylinder $K_{\vec{r}}(y)\times (s, t)\in \Omega_T$ and every piecewise smooth, cutoff function $\zeta(x)$ vanishing on $\partial K_{\vec{r}}(y)$ and such that $0\leqslant \zeta(x) \leqslant 1$ there holds
\begin{multline}\label{eq2.2}
\sum\limits_{i=1}^N\int\limits_s^t\int\limits_{K_{\vec{r}}(y)}\Big(\frac{\tau-s}{t-s}\Big)^\alpha
\frac{u^{(m_i-m_-)(p_i-1)}}{(u^{m_-}+\varepsilon^{m^-})^\beta}|D_i u^{m_-}|^{p_i}\,\zeta^{p_+}\,dx d\tau \\
\leqslant
\gamma \int\limits_{K_{\vec{r}}(y)\times\{t\}}(u+\varepsilon)^{1-m_-(\beta-1)}\,dx+
\gamma \sum\limits_{i=1}^N\int\limits_s^t\int\limits_{K_{\vec{r}}(y)}(u+\varepsilon)^{\lambda_i-m_-(\beta-1)}\,|D_i \zeta|^{p_i}\,dx d\tau.
\end{multline}
\end{lemma}
\begin{proof}
Test \eqref{eq1.3} by $(u^{m_-}+\varepsilon^{m_-})^{1-\beta}\big(\frac{\tau-s}{t-s}\big)^\alpha \zeta^{p_+}$ and integrate over $K_{\vec{r}}(y)\times(s, t)$. Using conditions \eqref{eq1.2} and the Young inequality we obtain
\[K_1(\beta-1)\sum\limits_{i=1}^N\int\limits_s^t\int\limits_{K_{\vec{r}}(y)}\Big(\frac{\tau-s}{t-s}\Big)^\alpha
\frac{u^{(m_i-m_-)(p_i-1)}}{(u^{m_-}+\varepsilon^{m_-})^\beta}|D_i u^{m_-}|^{p_i}\,\zeta^{p_+}\,dx d\tau \]
\begin{eqnarray*}
& \leqslant &
\gamma \int\limits_{K_{\vec{r}}(y)\times\{t\}}\int\limits_0^u\frac{ds}{(s^{m_-}+\varepsilon^{m_-})^{\beta-1}}\,\zeta^{p_+}\,dx\\
& & +
\gamma\sum\limits_{i=1}^N\int\limits_s^t\int\limits_{K_{\vec{r}}(y)}\Big(\frac{\tau-s}{t-s}\Big)^\alpha u^{(m_i-m_-)(p_i-1)}(u^{m_-}+\varepsilon^{m_-})^{p_i-\beta}|D_i\zeta|^{p_i}\,dx d\tau.
 \end{eqnarray*}
 From this, recalling the range of $\beta$ and using the inequalities
 \begin{eqnarray*}
 \gamma \int\limits_{K_{\vec{r}}(y)\times\{t\}}\int\limits_0^u\frac{ds}{(s^{m_-}+\varepsilon^{m_-})^{\beta-1}}\,\zeta^{p_+}\,dx & \leqslant & 
 \gamma \int\limits_{K_{\vec{r}}(y)\times\{t\}}\int\limits_0^u\frac{ds}{(s+\varepsilon)^{m_-(\beta-1)}}\,dx 
 \\
 & \leqslant &\gamma(\beta)\int\limits_{K_{\vec{r}}(y)}(u+\varepsilon)^{1-m_-(\beta-1)}\,dx,
 \end{eqnarray*}
 and
 \[
\sum\limits_{i=1}^N\int\limits_s^t\int\limits_{K_{\vec{r}}(y)}\Big(\frac{\tau-s}{t-s}\Big)^\alpha u^{(m_i-m_-)(p_i-1)}(u^{m_-}+\varepsilon^{m_-})^{p_i-\beta}|D_i\zeta|^{p_i}\,dx d\tau \]
\begin{eqnarray*}
& \leqslant& 
\sum\limits_{i=1}^N\int\limits_s^t\int\limits_{K_{\vec{r}}(y)}(u^{m_-}+\varepsilon^{m_-})^{\frac{m_i}{m_-}(p_i-1)+1-\beta}|D_i \zeta|^{p_i}\,dx d\tau \\
& \leqslant &   \gamma \, \sum\limits_{i=1}^N\int\limits_s^t\int\limits_{K_{\vec{r}}(y)}(u+\varepsilon)^{\lambda_i-m_-(\beta-1)}|D_i\zeta|^{p_i}\,dx d\tau,
 \end{eqnarray*}
 we arrive at the required \eqref{eq2.2}.
\end{proof}

\begin{lemma}\label{lem2.5}
Let $u$ be a non-negative, locally bounded, local weak sub-solution to \eqref{eq1.1}-\eqref{eq1.2} in $\Omega_T$. Construct the cylinder $K_{\vec{r}}(y)\times (s, t)\subset \Omega_T$ and assume that
\begin{equation}\label{eq2.3}
u\leqslant \mu^+,\quad \text{in}\quad K_{\vec{r}}(y)\times (s, t) \ .
\end{equation}
Let $\alpha>0$ and take $1<\beta <2 $. Then there exists $\gamma >0$, depending only on the data and $\beta$, such that for every $\varepsilon>0$ and every piecewise smooth cutoff function $\zeta(x)$, vanishing on $\partial K_{\vec{r}}(y)$ and such that $0\leqslant \zeta(x) \leqslant 1$, there holds
\[
\sum\limits_{i=1}^N\int\limits_s^t\int\limits_{K_{\vec{r}}(y)}\Big(\frac{\tau-s}{t-s}\Big)^\alpha
\frac{u^{(m_i-m_-)(p_i-1)}}{((\mu^+)^{m_-}-u^{m_-} +\varepsilon^{m_-})^\beta}|D_i u^{m_-}|^{p_i}\,\zeta^{p_+}\,dx d\tau\]
\begin{eqnarray}\label{eq2.4}
& \leqslant &
\gamma \, (\mu^+)^{1- m_-}\, \int\limits_{K_{\vec{r}}(y)
\times\{t\}} ((\mu^+)^{m_-}-u^{m_-} +\varepsilon^{m_-})^{2-\beta} \,\zeta^{p_+}\, dx \nonumber \\
&  & + \gamma \sum\limits_{i=1}^N[\mu^+]^{\lambda_i-m_-(p_i-1)}\int\limits_s^t\int\limits_{K_{\vec{r}}(y)}((\mu^+)^{m_-}-u^{m_-} +\varepsilon^{m_-})^{p_i-\beta}\,|D_i \zeta|^{p_i}\,dx d\tau.
\end{eqnarray}

\end{lemma}

\begin{proof}
 Test \eqref{eq1.3} by $((\mu^+)^{m_-}-u^{m_-} +\varepsilon^{m_-})^{1-\beta}\big(\frac{\tau-s}{t-s}\big)^\alpha \zeta^{p_+}$ and integrate over $K_{\vec{r}}(y)\times(s, t)$. Using the structure conditions \eqref{eq1.2} and the Young inequality we obtain
\[ 
\sum\limits_{i=1}^N\int\limits_s^t\int\limits_{K_{\vec{r}}(y)}\Big(\frac{\tau-s}{t-s}\Big)^\alpha
\frac{u^{(m_i-m_-)(p_i-1)}}{((\mu^+)^{m_-}-u^{m_-} +\varepsilon^{m_-})^\beta}|D_i u^{m_-}|^{p_i}\,\zeta^{p_+}\,dx d\tau \]
 \begin{eqnarray*}
  & \leqslant & \int\limits_{K_{\vec{r}}(y)
\times\{t\}} \left(\int_{u}^{\mu^+} ((\mu^+)^{m_-}-s^{m_-} +\varepsilon^{m_-})^{1-\beta} \ ds\right) \,\zeta^{p_+}\,dx \\
& & + \gamma \sum\limits_{i=1}^N \int\limits_s^t\int\limits_{K_{\vec{r}}(y)}u^{(m_i-m_-)(p_i-1)}((\mu^+)^{m_-}-u^{m_-} +\varepsilon^{m_-})^{p_i-\beta}\,|D_i \zeta|^{p_i}\,dx d\tau \\
 & \leqslant & \gamma \, (\mu^+)^{1- m_-}\, \int\limits_{K_{\vec{r}}(y)
\times\{t\}} ((\mu^+)^{m_-}-u^{m_-} +\varepsilon^{m_-})^{2-\beta} \,\zeta^{p_+}\, dx \\
& & +\gamma \sum\limits_{i=1}^N  (\mu^{+})^{(m_i-m_-)(p_i-1)}\int\limits_s^t\int\limits_{K_{\vec{r}}(y)}((\mu^+)^{m_-}-u^{m_-} +\varepsilon^{m_-})^{p_i-\beta}\,|D_i \zeta|^{p_i}\,dx d\tau  \ ,
\end{eqnarray*}
since, for $1<\beta<2$, 
\begin{eqnarray*}
    \int_{u}^{\mu^+} ((\mu^+)^{m_-}-s^{m_-} +\varepsilon^{m_-})^{1-\beta} \ ds 
    & \leqslant & \int_{\frac{u+\mu^+}{2}}^{\mu^+} ((\mu^+)^{m_-}-s^{m_-} +\varepsilon^{m_-})^{1-\beta} \ ds  \\
    & = & \frac{1}{m_-} 
    \int_{(\frac{u+\mu^+}{2})^{m_-}}^{(\mu^+)^{m_-}} ((\mu^+)^{m_-}-r +\varepsilon^{m_-})^{1-\beta} \ r^{\frac{1-m_-}{m_-}} \ dr\\
    & \leqslant &  \frac{\gamma(m_-)}{2-\beta} \, (\mu^+)^{1- m_-}\, ((\mu^+)^{m_-}-u^{m_-} +\varepsilon^{m_-})^{2-\beta} \ .
\end{eqnarray*}

\end{proof}

\section{Integral Harnack-type estimates}\label{S:harnackineq}

\noindent The principal focus of this section is the presentation (and proof) of certain inequalities known in the literature as integral (or $L^1-L^1$) Harnack-type estimates; we will do it separately for the singular case $\lambda_+< 1$ and for $p_+<2$, without additional assumptions on the gap between $\lambda_+$ and $\lambda_-$.

\begin{propo}\label{pr6.1}
Let $u$ be a non-negative, local weak solution to \eqref{eq1.1}- \eqref{eq1.2} and assume that $\lambda_{+}<1$. Then
there exists $\gamma >0$, depending only on the data, such that
\begin{equation}\label{eq6.1}
\sup\limits_{s<\tau<t}\fint\limits_{K_{\vec{r}}(y)} u(x, \tau)\,dx\leqslant \gamma\inf\limits_{s<\tau<t}\fint\limits_{K_{2\vec{r}}(y)} u(x, \tau)\,dx+
\gamma \sum\limits_{i=1}^N \Big(\frac{t-s}{r^{p_i}_i}\Big)^{\frac{1}{1-\lambda_i}},
\end{equation}
provided that $K_{2\vec{r}}(y)\times (s, t) \subset \Omega_T.$
\end{propo}

\begin{proof} 
Assume that $s=0$. For $j=0,1,\ldots$ set
\[
\vec{r}_{j}:= \sum\limits^{j}\limits_{l=0} \frac{\vec{r}}{2^{l}}, \quad
\bar{\vec{r}}_{j}:= \frac{r_{j}+r_{j+1}}{2}, \quad K_{j}:= K_{\vec{r}_{j}}(y), \quad \bar{K_{j}}:= K_{\bar{\vec{r}}_{j}}(y)
\]
and let $\zeta_j(x)\in C^1_0(\bar{K}_j)$, $\zeta_j(x)=1$ in $K_j$, $0\leqslant \zeta_j(x)\leqslant 1$, $|D_i \zeta_{j}|\leqslant \gamma \ \dfrac{2^{j}}{r_i}$, $i=1,...,N$. Fix $t_1,t_2\in(0,t)$.  For $t_1<t_2$,
testing \eqref{eq1.3} with $\zeta_j(x)$ and integrating over $\bar{K}_j \times (t_{1}, t_{2})\subset \Omega \times (0, t)$, we obtain 
\begin{eqnarray*}
\fint\limits_{\bar{K_{j}}} u(x, t_{2})\,\zeta_{j}\, dx
& \leqslant & 
\fint\limits_{\bar{K_{j}}} u(x, t_{1})\,\zeta_{j}\, dx\\
& & 
+ \gamma \sum\limits_{i=1}^N \frac{2^{j}}{r_i} \int\limits_{t_{1}}\limits^{t_{2}}
\fint\limits_{\bar{K_{j}}}  u^{\frac{(m_i-m_-)(p_i-1)}{p_i}} \left(\sum_{l=1}^N u^{(m_l-m_-)(p_l-1)}|D_l u^{m^-}|^{p_l} \right)^{\frac{p_i-1}{p_i}}\ dx d\tau.
\end{eqnarray*}
For $t_1>t_2$, we apply the weak formulation with the test function $-\zeta_j(x)$. Consequently, the inequality holds for arbitrary
$t_1,t_2\in(0,t)$.
By choosing $t_{1}$ such that 
\[
\fint\limits_{K_{2\vec{r}}(y)} u(x, t_{1})\,dx= \inf\limits_{0<\tau<t} \fint\limits_{K_{2\vec{r}}(y)\times\{\tau\}} u(x,\tau)\,dx \ ,\]
observing that
\[
\fint\limits_{\bar{K}_j}u(x,t_1)\,\zeta_j\,dx
\leqslant \gamma
\inf\limits_{0<\tau<t} \fint\limits_{K_{2\vec{r}}(y)\times\{\tau\}} u(x,\tau)\,dx,  \quad \bar{K}_j\subset K_{2\vec{r}}(y) \ ,
\]
and using the arbitrariness of $t_2$ we arrive at 
\begin{eqnarray*}
    J_j & \leqslant & \gamma 
\inf\limits_{0<\tau<t} \fint\limits_{K_{2\vec{r}}(y)\times\{\tau\}}u(x,\tau)\, dx \\
& & + \gamma \, 2^j \sum\limits_{i=1}^N \frac{1}{r_i} \int\limits_{t_{1}}\limits^{t_{2}}
\fint\limits_{\bar{K_{j}}}  u^{\frac{(m_i-m_-)(p_i-1)}{p_i}} \left(\sum_{l=1}^N u^{(m_l-m_-)(p_l-1)}|D_l u^{m^-}|^{p_l} \right)^{\frac{p_i-1}{p_i}}\ dx d\tau
\end{eqnarray*} 
being $\displaystyle{J_j:=\sup\limits\limits_{0<\tau<t}\,\,\fint\limits_{K_j}u(x,\tau)\,dx}$.
Our next goal is to derive an estimate from above for the second term of the right hand side of the previous inequality. For that purpose, we start by applying H\"older's inequality
\[\sum\limits_{i=1}^N \int\limits_{t_{1}}\limits^{t_{2}}
\fint\limits_{\bar{K_{j}}}  u^{\frac{(m_i-m^-)(p_i-1)}{p_i}} \left(\sum_{l=1}^N u^{(m_l-m_-)(p_l-1)}|D_l u^{m^-}|^{p_l} \right)^{\frac{p_i-1}{p_i}}\ dx d\tau \leq \sum\limits_{i=1}^N \frac{1}{r_i}  I_1^{\frac{p_i-1}{p_i}} \, I_2^{\frac{1}{p_i}}, \]
where
\[ I_1= \sum_{l=1}^N \int\limits_{0}\limits^{t}
\fint\limits_{\bar{K_{j}}}  \left(\frac{\tau}{t}\right)^\alpha \, \frac{u^{(m_l-m_-)(p_l-1)}}{ (u^{m_-}+\epsilon^{m_-})^\beta}  |D_l u^{m^-}|^{p_l} \zeta^{p_+} \ dx d\tau \ , \]
for $\zeta(x)\in C^1_0(K_{j+1})$, $\zeta(x)=1$ in $\bar{K}_j$, $0\leqslant \zeta_j(x)\leqslant 1$, $|D_i \zeta|\leqslant \gamma \ \dfrac{2^{j}}{r_i}$, $i=1,...,N$, and
\[I_2= \int\limits_{0}\limits^{t}
\fint\limits_{\bar{K_{j}}}  \left(\frac{\tau}{t}\right)^{-\alpha(p_i-1)} u^{(m_i-m_-)(p_i-1)} (u^{m_-}+\epsilon^{m_-})^{\beta(p_i-1)} \, dx d\tau \ .\]
By using Lemma \ref{lem2.4}, for $\displaystyle{\epsilon= \sum_{i=1}^N \left(\frac{t}{r_i^{p_i}}\right)^{\frac{1}{1-\lambda_i}}}$ and $1<\beta<1 +\frac{1}{m_-}$, we get
\begin{eqnarray}
I_1 & \leqslant & \gamma  \fint\limits_{K_{j+1}\times\{t\}}(u+\varepsilon)^{1-m_-(\beta-1)}\,dx+
 \gamma  \sum\limits_{l=1}^N\int\limits_s^t\fint\limits_{K_{j+1}}(u+\varepsilon)^{\lambda_l-m_-(\beta-1)}\,|D_l \zeta|^{p_l}\,dx dt  \nonumber \\
& \leqslant & \gamma \sup\limits_{0<\tau<t}
\fint\limits_{K_{j+1}\times\{\tau\}}(u+\varepsilon)^{1-m_-(\beta-1)}\,dx \left[ 1+ \sum\limits_{l=1}^N \frac{\epsilon^{\lambda_l-1}t}{r_i^{p_l}}\right] \nonumber \\
& \leqslant & \gamma  \sup\limits_{0<\tau<t}
\fint\limits_{K_{j+1}\times\{\tau\}}(u+\varepsilon)^{1-m_-(\beta-1)}\,dx \leqslant
\gamma  \left(\sup\limits_{0<\tau<t}
\fint\limits_{K_{j+1}\times\{\tau\}}u+\varepsilon \, dx\right)^{1-m_-(\beta-1)} \ . 
\end{eqnarray}
As for $I_2$, considering $1<\beta<1 +\frac{1-\lambda_+}{m_-(p_+-1)}$ and taking $\alpha <\frac{1}{p_+-1} $ we obtain
\begin{eqnarray}
I_2 & \leqslant & \gamma  \int\limits_{0}\limits^{t}
\fint\limits_{\bar{K_{j+1}}}  \left(\frac{\tau}{t}\right)^{-\alpha(p_i-1)} (u+\epsilon)^{(m_i+m_-(\beta-1))(p_i-1)} \, dx d\tau \nonumber \\
& \leqslant & \gamma\,  t\, \sup\limits_{0<\tau<t}
\fint\limits_{K_{j+1}\times\{\tau\}}(u+\epsilon)^{(m_i+m_-(\beta-1))(p_i-1)} \, dx , \nonumber \\
& \leqslant & \gamma\,  t\, \left(\sup\limits_{0<\tau<t}
\fint\limits_{K_{j+1}\times\{\tau\}} u+\epsilon \, dx \right)^{(m_i+m_-(\beta-1))(p_i-1)} \ .
\end{eqnarray}
Gathering these estimates altogether and using Young's inequality with $\delta\in (0,1)$, we arrive at 
\begin{eqnarray*}
J_j 
& \leqslant & \gamma 
\inf\limits_{0<\tau<t} \fint\limits_{K_{2\vec{r}}(y)\times\{\tau\}}u(x,\tau)\, dx  + \gamma \, 2^j  
    \sum\limits_{i=1}^N \frac{t^{\frac{1}{p_i}}}{r_i} \left(\sup\limits_{0<\tau<t}
\fint\limits_{K_{j+1}\times\{\tau\}}u+\varepsilon \, dx\right)^{\frac{(m_i+1)(p_i-1)}{p_i}} \\
& \leqslant & \gamma 
\inf\limits_{0<\tau<t} \fint\limits_{K_{2\vec{r}}(y)\times\{\tau\}}u(x,\tau)\, dx + \delta \sup\limits_{0<\tau<t}
\fint\limits_{K_{j+1}\times\{\tau\}}u(x,\tau) \, dx + 
\gamma 2^{\gamma j} \delta^{-\gamma} \epsilon
\end{eqnarray*}
and therefore, by iteration and choosing $\delta$ so small that $2^\gamma\delta<1$, we arrive at 
\begin{eqnarray*}
\sup\limits_{0<\tau<t}\,\,\fint\limits_{K_{\vec{r}}(y)} u(x,\tau)\,dx & \leqslant & \delta^j J_j + \gamma \delta^{-\gamma} \Big(\inf\limits_{0< \tau <t} \fint\limits_{K_{2\vec{r}}(y)\times\{\tau\}} u\ dx
+ \varepsilon\Big)\sum\limits_{i=0}^{j-1} (\delta2^{\gamma})^i\\
& \leqslant & \delta^j \gamma\sup\limits_{0<\tau<t}\fint\limits_{K_{2\vec{r}}(y)}u(x,\tau)\, dx  + \gamma \Big(\inf\limits_{0< \tau <t} \fint\limits_{K_{2\vec{r}}(y)\times\{\tau\}} u\ dx
+ \sum_{i=1}^N \left(\frac{t}{r_i^{p_i}}\right)^{\frac{1}{1-\lambda_i}} \Big) \ . 
\end{eqnarray*}
The proof is complete once we pass to the limit as $j\to\infty$. 

\end{proof}

\begin{propo}\label{pr6.2}
Let $u$ be a non-negative, locally bounded, local weak solution to \eqref{eq1.1}-\eqref{eq1.2} and $p_+<2$. Construct the cylinder
$K_{2\vec{r}}(y)\times (s, t) \subset \Omega_T$ and let 
$$u\leqslant \mu^+ \quad \text{in}\quad K_{2\vec{r}}(y)\times (s, t).$$
Then there exists $\gamma >0$, depending only on the data, such that 
\begin{equation}\label{eq6.2}
\sup\limits_{s\leqslant \tau\leqslant t}\fint\limits_{K_{\vec{r}}(y)}\mu^+-u(x, \tau) \, dx \leqslant \gamma
\inf\limits_{s<\tau <t}\fint\limits_{K_{2\vec{r}}(y)} \mu^+-u(x, \tau) \,dx +\gamma \sum\limits_{i=1}^N
\Big((\mu^+)^{\lambda_i-p_i+1} \, \frac{t-s}{r^{p_i}_i}\Big)^{\frac{1}{2-p_i}}.
\end{equation}
\end{propo}

\noindent The proof is quite similar to the one of Proposition \ref{pr6.1}, nevertheless and for the sake of clarify we decided to present it.

\begin{proof} As in the previous proof, consider $s=0$ and construct the following sequences of radii, cubes and cutoff functions, $ j=0,1,\cdots$,
\[
\vec{r}_{j}:= \sum\limits^{j}\limits_{l=0} \frac{\vec{r}}{2^{l}}, \quad
\bar{\vec{r}}_{j}:= \frac{r_{j}+r_{j+1}}{2}, \quad K_{j}:= K_{\vec{r}_{j}}(y), \quad \bar{K_{j}}:= K_{\bar{\vec{r}}_{j}}(y)
\]
\[ \zeta_j(x)\in C^1_0(\bar{K}_j) \ , \ \ \zeta_j(x)=1 \ \mbox{in} \  K_j \ , \ \ 0\leqslant \zeta_j(x)\leqslant 1 \ , \ \ |D_i \zeta_{j}|\leqslant \gamma \ \dfrac{2^{j}}{r_i}, i=1,...,N. \] 

\noindent  Testing \eqref{eq1.3} with $\varphi=\zeta_j(x)$
and integrating over $\bar K_j\times(t_1,t_2)$,
we rewrite the time term as
\[\int\limits_{\bar{K}_j} u(x,t)\,\zeta_j(x)\,dx \Big|_{t_1}^{t_2}=
\int\limits_{\bar{K}_j} (\mu^+ - u(x,t_1)) \zeta _j\ dx  -\int\limits_{\bar{K}_j}
(\mu^+ - u(x,t_2))
\,\zeta_j(x)\,dx\]
Thus, we obtain
\begin{multline*}
\fint_{\bar K_j}(\mu^+-u(x,t_2))
\zeta_j(x)\,dx
\leqslant
\fint_{\bar K_j}(\mu^+-u(x,t_1))
\zeta_j(x)\,dx\\
+\gamma 2^j\sum_{i=1}^N\frac{1}{r_i}
\int_{t_1}^{t_2}
\fint_{\bar K_j}
u^{\frac{(m_i-m_-)(p_i-1)}{p_i}}
\left(
\sum_{l=1}^N
u^{(m_l-m_-)(p_l-1)}
|D_lu^{m_-}|^{p_l}
\right)^{\frac{p_i-1}{p_i}}
\,dx\,d\tau
\end{multline*}
for arbitrary $t_1, t_2 \in (0, t)$.
More precisely, to obtain this estimate for times
$t_1 > t_2$, we use \eqref{eq1.3} with the test function $-\zeta_{j}(x)$.

\noindent Choose $t_1$ is such that 
\[
\fint\limits_{K_{2\vec{r}}(y)} \mu^+-u(x, t_{1})\,dx= \inf\limits_{0<\tau<t} \fint\limits_{K_{2\vec{r}}(y)\times\{\tau\}} \mu^+-u \, dx \ .\]
Then, due to the arbitrariness of $t_2$ and for 
$\displaystyle{J_j:=\sup\limits\limits_{0<\tau<t}\,\,\fint\limits_{K_j \times\{\tau\}}\mu^+ -u \,dx}$, we arrive at
\begin{eqnarray*}
    J_j & \leqslant & \gamma 
\inf\limits_{0<\tau<t} \fint\limits_{K_{2\vec{r}}(y)\times\{\tau\}}\mu^+-u \, dx \\
& & + \gamma \, 2^j \sum\limits_{i=1}^N \frac{1}{r_i} \int\limits_{t_{1}}\limits^{t_{2}}
\fint\limits_{\bar{K_{j}}}  u^{\frac{(m_i-m_-)(p_i-1)}{p_i}} \left(\sum_{l=1}^N u^{(m_l-m_-)(p_l-1)}|D_l u^{m^-}|^{p_l} \right)^{\frac{p_i-1}{p_i}}\ dx d\tau\\
& \leqslant & \gamma \inf\limits_{0<\tau<t} \fint\limits_{K_{2\vec{r}}(y)\times\{\tau\}}\mu^+-u \, dx + \gamma \sum\limits_{i=1}^N \frac{2^j}{r_i}  I_1^{\frac{p_i-1}{p_i}} \, I_2^{\frac{1}{p_i}},
\end{eqnarray*} 
where
\[ I_1= \sum_{l=1}^N \int\limits_{0}\limits^{t}
\fint\limits_{\bar{K_{j}}}  \left(\frac{\tau}{t}\right)^\alpha \, \frac{u^{(m_l-m_-)(p_l-1)}}{ ((\mu^+)^{m_-} - u^{m_-}+\epsilon^{m_-})^\beta}  |D_l u^{m^-}|^{p_l} \zeta^{p_+} \ dx d\tau \ , \]
\[I_2= \int\limits_{0}\limits^{t}
\fint\limits_{\bar{K_{j}}}  \left(\frac{\tau}{t}\right)^{-\alpha(p_i-1)} u^{(m_i-m_-)(p_i-1)} ((\mu^+)^{m_-} - u^{m_-}+\epsilon^{m_-})^{\beta(p_i-1)} \, dx d\tau \ \]
and $\zeta(x)\in C^1_0(K_{j+1})$, is a time-independent cutoff function verifying \[ 0\leqslant \zeta_j(x)\leqslant 1 \ , \quad  \zeta(x)=1, \ \mbox{in} \  \bar{K}_j \ , \quad |D_i \zeta|\leqslant \gamma \ \dfrac{2^{j}}{r_i} \ , i=1,...,N  \ .\]
By using Lemma \ref{lem2.5}, for $\displaystyle{\epsilon^{m_-}= \sum_{i=1}^N \left(\frac{(\mu^+)^{\lambda_i-1+m_-(2-p_i)}\, t}{r_i^{p_i}}\right)^{\frac{1}{2-p_i}}}$ and $1<\beta<2$, we get

\[ I_1 \leqslant \gamma \, 2^{\gamma j}\, (\mu^+)^{1-m_-} \,   \left(\sup\limits_{0<\tau<t}
\fint\limits_{K_{j+1}\times\{\tau\}}(\mu^+)^{m_-} - u^{m_-}+\epsilon^{m_-} \, dx\right)^{2-\beta} \]
$$ + \gamma \sum\limits_{l=1}^N[\mu^+]^{\lambda_l-m_-(p_l-1)}\int\limits_s^t\int\limits_{K_{\vec{r}}(y)}((\mu^+)^{m_-}-u^{m_-} +\varepsilon^{m_-})^{p_l-\beta}\,|D_l \zeta|^{p_l}\,dx d\tau.$$
and by considering $1<\beta<\frac{1}{p_+-1}$ and $\alpha< \frac{1}{p_+-1}$ we get
\[ 
I_2  \leqslant  \gamma\, (\mu^+)^{\lambda_i- m_-(p_i-1)}\,  t\, \left(\sup\limits_{0<\tau<t}
\fint\limits_{K_{j+1}\times\{\tau\}} (\mu^+)^{m_-} - u^{m_-}+\epsilon^{m_-} \, dx \right)^{\beta(p_i-1)} \ . \]
Gathering these estimates altogether and using Young's inequality with $\delta\in (0,1)$, we arrive at 
\begin{eqnarray*}
J_j 
& \leqslant & \gamma 
\inf\limits_{0<\tau<t} \fint\limits_{K_{2\vec{r}}(y)\times\{\tau\}}\mu^+-u \, dx \\
& & + \gamma \, 2^j  (\mu^+)^{1-m_-} 
    \sum\limits_{i=1}^N \frac{t^{\frac{1}{p_i}}}{r_i} \, (\mu^+)^{\frac{\lambda_i-1+m_-(p_i-1)}{p_i}}\, \left(\sup\limits_{0<\tau<t}
\fint\limits_{K_{j+1}\times\{\tau\}} ((\mu^+)^{m_-} - u^{m_-}+\epsilon^{m_-})\, dx \right)^{\frac{2(p_i-1)}{p_i}} \\
& \leqslant & \gamma 
\inf\limits_{0<\tau<t} \fint\limits_{K_{2\vec{r}}(y)\times\{\tau\}}\mu^+-u \, dx  + \delta (\mu^+)^{1-m_-} \sup\limits_{0<\tau<t}
\fint\limits_{K_{j+1}\times\{\tau\}} (\mu^+)^{m_-} - u^{m_-}\, dx \\
& & +  
\gamma 2^{\gamma j} \delta^{-\gamma} \epsilon^{m_-}(\mu^{+})^{1-m_-} \\
& \leqslant & \gamma 
\inf\limits_{0<\tau<t} \fint\limits_{K_{2\vec{r}}(y)\times\{\tau\}}\mu^+-u \, dx  + \delta \sup\limits_{0<\tau<t}
\fint\limits_{K_{j+1}\times\{\tau\}} \mu^+ - u \, dx +  
\gamma 2^{\gamma j} \delta^{-\gamma} \epsilon^{m_-}(\mu^{+})^{1-m_-} \ .
\end{eqnarray*}
Once we iterate and  choose $\delta$ so small that $2^\gamma\delta<1$
\begin{eqnarray*}
\sup\limits_{0<\tau<t}\,\,\fint\limits_{K_{\vec{r}}(y)\times\{\tau\}} \mu^+-u \,dx & \leqslant & \delta^j J_j \\
& & + \gamma \delta^{-\gamma} \Big(\inf\limits_{0< \tau <t} \fint\limits_{K_{2\vec{r}}(y)\times\{\tau\}} \mu^+-u\ dx
+ \epsilon^{m_-}(\mu^{+})^{1-m_-}\Big)\sum\limits_{i=0}^{j-1} (\delta2^{\gamma})^i\\
& \leqslant & \delta^j \gamma\sup\limits_{0<\tau<t}\fint\limits_{K_{2\vec{r}}(y)}\mu^+- u(x,\tau)\, dx  \\
& & + \gamma \Big(\inf\limits_{0< \tau <t} \fint\limits_{K_{2\vec{r}}(y)\times\{\tau\}} \mu^+-u\ dx + \sum_{i=1}^N \left(\frac{t \, (\mu^+)^{\lambda_i+1-p_i}}{r_i^{p_i}}\right)^{\frac{1}{2-p_i}} \Big) \ . 
\end{eqnarray*}
The desired inequality is obtained by letting $j \rightarrow +\infty$.
\end{proof}

\section{Proving Theorem \ref{th1.1}}\label{S:provingthe1.1}

\noindent Fix $(x_0, t_0)\in \Omega_T$, and consider $0<R<1$ such that 
\[Q_{8R}:= \{|x_i-x_{i,0}|<(8R)^{\frac{1}{2}},\,\,i=1, ..., N\}\times(t_0-8R, t_0)\subset \Omega_T.\]
If for all $0<\rho\leqslant R$
\[ \osc\limits_{Q_\rho} u\leqslant \,\rho^{\frac{p_+}{2}},\]
then there is nothing to prove. Assume that there exists $\rho_0 \in (0, R]$ such that
\begin{equation}\label{eq3.10}
\osc\limits_{Q_{\rho_0}} u\geqslant \rho_0^{\frac{p_+}{2}}.
\end{equation}

\subsection{Reducing the oscillation - general setting}

\noindent Fix $(x_0,t_0) \in \Omega_T$ and let $0<\rho<R$ and $\mu^+$, $\mu^-$ and $\omega$ be non-negative real numbers verifying
\[\mu^+\geqslant \sup\limits_{Q ^{\omega}_{\rho, \theta(\omega,\rho) }}u,\quad \mu^-:=\inf\limits_{Q ^{\omega}_{\rho,\theta(\omega,\rho) }} u ,\quad \mu^+-\mu^-\leqslant \omega,\]
where \[ Q ^{\omega}_{\rho, \theta(\omega,\rho)}= Q ^{\omega}_{\rho, \theta(\omega, \rho)}(x_0,t_0) \ , \quad \mbox{for} \quad  \  \theta(\omega, \rho) = b \, \omega^{1-\lambda_+} \rho^{p_+}  , \ \ b>0\ .\]
\begin{propo}\label{propgeneralset}
 Let $u$ be a non-negative, local weal solution to \eqref{eq1.1}-\eqref{eq1.2}. Consider $\lambda_+<1$  and $\varkappa >0$. Then there exists $b, \eta, \epsilon, \sigma_0 \in (0,1)$, depending only on the data, such that if
 \begin{equation}\label{condmu+}
\mu^+ \leqslant (1+\eta) \omega
\end{equation}
then 
\begin{equation}\label{redosc}
\osc\limits_{Q ^{\sigma_0 \omega}_{\epsilon \rho,  \theta(\sigma_0 \omega, \epsilon \rho)}} u \leqslant \sigma_0  \ \omega \ ,
\end{equation}
provided $\lambda_+-\lambda_- \leqslant \epsilon_*$, for some $\epsilon_*= \epsilon_*(data) \in (0,1)$.
\end{propo}

\begin{remark}
    We assume, without loss of generality, that $\mu^+ \geqslant \frac{\omega}{2}$. In fact, if not then \eqref{redosc} comes immediately for $\sigma_0=1/2$.
\end{remark}

\begin{proof}  The proof of this result is based on the study of the  following two alternatives: there exists $\nu \in (0,1)$, depending only upon the data and $b$ (to de defined), such that either
\begin{equation}\label{eq3.16}
\left| Q^{\omega}_{\frac{\rho}{2}, \theta(\omega, \rho/2)} \cap \{u\geqslant \mu^+-\frac{\omega}{4} \}\right| \leqslant \nu \left|Q^{\omega}_{\frac{\rho}{2}, \theta(\omega, \rho/2)}\right| \ ,
\end{equation}
or
\begin{equation}\label{eq3.17}
\left| Q^{\omega}_{\frac{\rho}{2},\theta(\omega, \rho/2)} \cap \{u\geqslant \mu^+-\frac{\omega}{4} \}\right| > \nu \left|Q^{\omega}_{\frac{\rho}{2}, \theta(\omega, \rho/2)}\right| \ .
\end{equation}

\paragraph{Analysis of the First Alternative.} Take $\nu=\frac{1}{\gamma} b^{\frac{N}{p}}$, by  De Giorgi type Lemma \ref{lem2.7} and Remark $2.2$,  we obtain from \eqref{eq3.16} that
$$u(x,t)\leqslant \mu^+-\frac{1}{8}\,\omega,\quad a.e. \, \, \text{in}\ Q^{\omega}_{\frac{\rho}{4},\theta(\omega, \rho/4)} $$
and hence
\begin{equation}\label{eq3.18}
\osc\limits_{Q^{\omega}_{\frac{\rho}{4},\theta(\omega, \rho/4) }} u \leqslant \frac{7}{8}\,\omega.
\end{equation}

\paragraph{Analysis of the Second Alternative.} Let $\nu=\frac{1}{\gamma} b^{\frac{N}{p}}$ and assume that \eqref{eq3.17} holds. Then, for some 
$\bar{t}\in (t_0- \theta(\omega, \rho/2), t_0)$
\begin{equation*}
\left| K^{\omega}_{\frac{\rho}{2}}(x_0)\cap\{u(\cdot, \bar{t})\geqslant \mu^+-\frac{\omega}{4}\}\right |\geqslant \nu \left|K^{\omega}_{\frac{\rho}{2}}(x_0)\right| \ .
\end{equation*}
We now use this measure theoretical information together with the bounds $ \frac{\omega}{2} \leqslant \mu^+ \leqslant 2 \omega $ and Proposition \ref{pr6.1} to get, for all $t_0-\theta(\omega,\rho) \leqslant t\leqslant t_0$,
\begin{eqnarray*}
\frac{b^{\frac{N}{p}}}{\tilde{\gamma}} 
\ \omega &\leqslant & \frac{\gamma}{4} b^{\frac{N}{p}}\omega  \leqslant  \nu \big(\mu^+-\frac{\omega}{4}\big)\leqslant
\fint\limits_{K^{\omega}_{\frac{\rho}{2}}(x_0)} u(x, \bar{t})\,dx  \\
& \leqslant & \tilde{\gamma}\inf\limits_{t_0-\theta(\omega, \rho/2) \leqslant t\leqslant t_0} \fint\limits_{K^{\omega}_{\rho}(x_0)} u(x, t)\,dx +\tilde{\gamma}\omega \ \sum\limits_{i=1}^N b^{\frac{1}{1-\lambda_i}} \\
& \leqslant & \tilde{\gamma} \epsilon \omega + \tilde{\gamma} \omega \, \frac{|K^{\omega}_{\rho}(x_0)| \cap \{ u(\cdot, t) \geqslant \epsilon \omega\} |}{|K^{\omega}_{\rho}(x_0)|} +\tilde{\gamma} \omega b^{\frac{N}{p}+\frac{\varkappa}{p(1-p|\frac{\lambda}{p}|)}}\sum\limits_{i=1}^N b^{\frac{1}{1-\lambda_i}-\frac{1}{1-p|\frac{\lambda}{p}|}} \\
 & \leqslant & \frac{b^{\frac{N}{p}}}{2 \tilde{\gamma}} 
\ \omega + \tilde{\gamma} \omega \, \frac{|K^{\omega}_{\rho}(x_0)| \cap \{ u(\cdot, t) \geqslant \epsilon \omega\} |}{|K^{\omega}_{\rho}(x_0)|} 
\end{eqnarray*}
by choosing
\[ b^{\frac{\varkappa}{p(1-p|\frac{\lambda}{p}|)}}\leqslant \frac{1}{8 \tilde{\gamma}^2} , \qquad  \epsilon= \frac{b^{\frac{N}{p}}}{4\tilde{\gamma}^2} \] 
and $\epsilon_{\star} $ so small that 
\[ \sum\limits_{i=1}^N b^{\frac{1}{1-\lambda_i}-\frac{1}{1-p|\frac{\lambda}{p}|}}\leqslant N b^{-\frac{\lambda_+-\lambda_-}{(1-\lambda_+)(1-p|\frac{\lambda}{p}|)}} \leqslant N b^{-\frac{\epsilon_*}{(1-\lambda_+)(1-p|\frac{\lambda}{p}|)}}\leqslant 2 \ . \]
Thereby, for all $t_0-\theta(\omega,\rho) \leqslant t\leqslant t_0$, 
\[ |K^{\epsilon \omega}_{\rho}(x_0)| \cap \{ u(\cdot, t) \geqslant \epsilon \omega\} | \geqslant \alpha |K^{\epsilon \omega}_{\rho}(x_0)| \ , \quad \alpha=\alpha(data) \in (0,1) \ , \] 
and finally, by applying  Lemma \ref{lem3.2}, we conclude
\begin{equation*}
u(x, t)\geqslant\sigma\,\omega,\quad (x, t)\in {K}^{\sigma \omega}_{\rho}(x_0) \times (t_0-\frac{1}{2}\theta(\omega, \rho)\, t_0), \qquad \mbox{for} \ \ \sigma= \dfrac{\epsilon}{C_\star},
\end{equation*}
Hence, recalling \eqref{condmu+},
\begin{equation}\label{eq3.21}
\osc\limits_{Q^{\sigma \omega}_{\rho, \frac{1}2 \theta(\omega, \rho)}} u\leqslant (1+\eta -\sigma )\omega=(1-\frac{1}{2}\sigma )\omega,
\end{equation}
provided that $\eta=\frac{1}{2}\sigma$. 

\paragraph{Combining the reduction of the oscillations.} 
Getting \eqref{eq3.18} and \eqref{eq3.21} altogether, we arrive at
\begin{equation}\label{eq3.22}
\osc\limits_{{Q}^{\sigma_0\omega}_{\epsilon_0 \rho, \theta (\epsilon_0 \rho, \sigma_0 \omega)}} u\leqslant \sigma_0 \,\omega \ .
\end{equation}
for 
\[ \epsilon_0= \frac{1}{4} \,\sigma_0^{\frac{\lambda_+-\lambda_-}{p_+}}  \qquad \mbox{and} \qquad 
\sigma_0= \max \left\{1- \frac{\sigma}{2} , \frac{7}{8} \right\}  =  1- \frac{\sigma}{2} >\frac{1}{4} > \sigma= \frac{\epsilon}{C_\star} \ . \]

\end{proof}

\noindent Proposition \ref{propgeneralset} will be the key ingredient to prove the main result leading to the H\"older continuity. That is precisely the content of the section to come.

\subsection{Main Proposition - iterative procedure}

\noindent Recall that we are working in the following geometric context: $(x_0, t_0)\in \Omega_T$ is a fixed point; $0<R<1$ is such that 
\[Q_{8R}:= \{|x_i-x_{i,0}|<(8R)^{\frac{1}{2}},\,\,i=1, ..., N\}\times(t_0-8R, t_0)\subset \Omega_T,\]
and $\rho_0 \in (0, R]$ is such that 
\begin{equation}\label{neweq3.10}
\osc\limits_{Q_{\rho_0}} u\geqslant \rho_0^{\frac{p_+}{2}} \ .
\end{equation}
Define the positive number
\begin{equation}\label{omega0}
\omega_0:= \max\left\{  \rho_0^{\frac{p_+}{2}} , \mu^{+}\right\} \qquad \mu^+:= \sup \limits_{Q_{\rho_0}} u  . 
\end{equation}

\begin{propo}\label{pr3.1}
There exist $\epsilon$, $b$, $\sigma_0 \in (0, 1)$ such that, for $j=0, 1, 2, ...$, 
\begin{equation}\label{eq3.15}
Q_{j+1} \subset Q_j \qquad \mbox{and}  \qquad  \mu^+_j-\mu^-_j\leqslant \omega_j \ , 
\end{equation}
for 
\[
\left\{
\begin{array}{l}
\mu^+_j:=\sup\limits_{Q_j} u , \quad  \mu^-_j:=\inf\limits_{Q_j} u \ ,\\
Q_j= Q^{\omega_j}_{\rho_j, \theta_j} \ , \quad \mbox{for} \quad \rho_j:=\epsilon^j\,\rho_0, \quad  \omega_j:=\sigma_0^j \ \omega_0 \ , \quad \theta_j= \theta(\omega_j, \rho_j)= b \ \rho_j^{p_+} \ \omega_j^{1-\lambda_+}\ ,
\end{array}
\right.
\]
provided $\lambda_+-\lambda_- \leqslant \epsilon_*$, for some $\epsilon_*= \epsilon_*(data) \in (0,1)$.
\end{propo}

\begin{proof} 
The proof is performed by induction on the natural variable $j$. Moreover, having a sequence of nested and shrinking cylinders $\left\{Q_j\right\}_j$ is a straightforward consequence of its construction once we take $\epsilon$ even small, if necessary, so that $\epsilon < \sigma_0^{\frac{\lambda_+-\lambda_-}{p_+}}$.

\paragraph{Step $j=0$.} According to \eqref{omega0} and to the fact that u in non-negative, we immediately get that $ \osc\limits_{Q_{\rho_0}} u \leqslant \omega_0$ thereby the desired reduction of the oscillation $\osc\limits_{Q_{0}} u \leqslant \omega_0$ is obtained by showing that $Q_0 \subset Q_{\rho_0}$. Observe that, due to \eqref{neweq3.10}, the closeness of $\lambda_i$, the smallness of $\rho_0$ and choosing $b$ even smaller (if needed), 
\[ \frac{{\rho_0}^{\frac{p_+}{p_i}}}{\omega_0^{\frac{\lambda_+-\lambda_i}{p_i}}}  \leqslant \rho_0^{\frac{p_+}{p_i} (1- \frac{\lambda_+-\lambda_i}{2})} \leq \rho_0^{\frac{1}{2}}\]
and
\[ b \omega_0^{1-\lambda_+} \rho_0^{p_+} \leqslant  \rho_0^{p_+}  \leqslant \rho_0 \]
This shows that $
Q_0 \subset Q_{\rho_0}$, hence 
\[\osc_{Q_0} u \leqslant \osc\limits_{Q_{\rho_0}} u \leqslant \omega_0 \]
and also 
\[\mu^{+}_0 = \sup_{Q_0} u \leqslant  \sup_{Q_{\rho_0}} u \leqslant \omega_0 < (1+\eta) \omega_0 \ . \] 
This is the starting point to go to {\textbf{Step} $j=1$}. In fact, since \eqref{condmu+} holds in the cylinder $Q_0$, we can make use of  Proposition \ref{propgeneralset} to determine $b, \epsilon, \sigma \in (0,1)$ such that
\[\osc_{Q_1} u  = \osc\limits_{{Q}^{\sigma_0\omega_0}_{\epsilon \rho_0, \theta (\epsilon \rho_0, \sigma_0 \omega_0)}} u \leqslant \sigma_0 \,\omega = \omega _1 \ , \qquad \mbox{for} \ \ \eta= \frac{\sigma}{2} \ ; \ \ \sigma_0=1-\frac{\sigma}{2} \ . \]
This procedure can be repeated as long as
\begin{equation}\label{eq3.23}
\mu^+_j\leqslant (1+\frac{\sigma}{2} )\,\omega_j \ ,
\end{equation}
leading to the reduction of the oscillation  \eqref{eq3.15}, over the cylinder $Q_j$. 

\noindent Assume now that, at some step $j_0$, \eqref{eq3.23} fails. Then 
\[ \mu^-_{j_0}\geqslant \frac{\sigma/2}{1+\sigma/2}\mu^+_{j_0} \ \]
and since \eqref{eq3.23} holds in the previous step $j_0-1$ 
\begin{equation}\label{eq3.25}
\frac{\sigma}{2} \,\, \omega_{j_0}\leqslant \frac{\sigma/2}{1+\sigma/2}\, \mu^+_{j_0}\leqslant \mu^-_{j_0}\leqslant \mu^+_{j_0-1}\leqslant \left(1+\frac{\sigma}{2}\right)\,\omega_{j_0-1}=\frac{1+\sigma/2}{\sigma_0}\, \omega_{j_0}.
\end{equation}
Consider the change of variables given by 
\[ z_i=(x_i-x_{i,0}) \ \omega^{\frac{\lambda_+-\lambda_i}{p_i}}_{j_0},\quad i=1, ..., N, \qquad \tau= t \ \omega_{j_0}^{\lambda_+-1} \]
which maps the cylinder $Q_{j_0}$ into  $Q_{\rho_{j_0}, b \rho_{j_0}^{p_+}}(0, \tau_0) $, being $\tau_0= t_0 \ \omega_{j_0}^{\lambda_+-1}$; and introduce the function $v(z, \tau)= \dfrac{u(x,t)}{\omega_{j_0}}$. Taking \eqref{eq3.25} into consideration, we observe that 
\[ \frac{\sigma}{2} \leqslant v(z,\tau) \leqslant 1+ \frac{\sigma}{2} \ , \quad (z,\tau) \in Q_{\rho_{j_0}, b \rho_{j_0}^{p_+}}(0, \tau_0) \]
and therefore, within the cylinder $Q_{\rho_{j_0}, b \rho_{j_0}^{p_+}}(0, \tau_0)$, $v$ satisfies the general anisotropic p-Laplacian 
\begin{equation}\label{eq3.6}
\begin{cases}
v_\tau-\sum\limits_{i=1}^N D_i \bar{a}_i(z, \tau, v, D v)=0,\quad (z,\tau)\in Q_{\rho_{j_0}, b \rho_{j_0}^{p_+}}(0, \tau_0) ,\\
\sum\limits_{i=1}^N \bar{a}_i(z, \tau, v, D v)\,D_i v \geqslant \bar{K}_1 \sum\limits_{i=1}^N \big|D_i v\big|^{p_i},\,\,p_- >1,\\
|\bar{a}_i(z,\tau, v, D v)|\leqslant \bar{K}_2  \,\big|D_i v\big|^{p_{i-1}},\quad i=1, ..., N,
\end{cases}
\end{equation}
for positive constants $\bar{K}_1$, $\bar{K}_2$, depending on the data and on $\sigma, \lambda_-, \lambda_+$. The study of the local continuity of $v$ can be found along \cite[section 7]{CiaHenSavSkr3} - as usual, the approach to the study of the  degenerate case,  $p_+>2$, differs from the one considered to the singular case, $1< p_+<2$.

\noindent We arrive at the required estimate \eqref{eq3.15} once we return to the original variables $(x,t)$ and function $u$.
\end{proof}

\noindent Now the H\"{o}lder continuity follows by standard iterative arguments, completing thereby the proof of Theorem \ref{th1.1}.

\section{Proving Theorem \ref{th1.2}} \label{S:provingthe1.2}

\noindent The beginning of the proof  is similar to that of Theorem \ref{th1.1}. In fact, we start by fixing $(x_0, t_0)\in \Omega_T$ and consider $0<R<1$ such that 
$$Q_{8R}:= \{|x_i-x_{i,0}|<(8R)^{\frac{1}{2}},\,\,i=1, ..., N\}\times(t_0-8R, t_0)\subset \Omega_T \ .$$
If for all $0<\rho\leqslant R$
$$\osc\limits_{Q_\rho} u\leqslant \mathcal{A}\,\rho^\beta,$$
with $\mathcal{A}> 1$ and $\beta>0$ to be fixed, then there is nothing to prove. If not, then there exists $\rho_0 \in (0, R]$ such that
\begin{equation}\label{eq4.5}
\osc\limits_{Q_{\rho_0}} u\geqslant \mathcal{A}\,\rho^{\beta}_0.
\end{equation}
Define $\displaystyle{\omega_0:= \max \left\{\mathcal{A}\,\rho^{\beta}_0, \sup_{Q_{\rho_0}} u \right\}}$ and construct the cylinder 
\begin{equation}\label{Qomega0}
Q^{\omega_0}_{\rho_0,\theta(\omega_0, \rho_0) }:= K^{\omega_0}_{\rho_0}(x_0)\times (t_0- \theta(\omega_0, \rho_0), t_0) \ , \quad \theta(\omega_0, \rho_0)= b \ \omega^{1-\lambda_+}_0 \rho^{p_+}_0, \ \ b>0 \ .
\end{equation}
 By considering 
\[ \mathcal{A} = \max \left\{1, b^{\frac{1}{\lambda_+-1}}\right\} \quad \mbox{ and } \quad \beta= \min\left\{ \frac{p_+}{2(\lambda_+-\lambda_-)}, \frac{p_+-1}{2(\lambda_+-1)} \right\}\]
we get
\[ Q^{\omega_0}_{\rho_0,\theta(\omega_0, \rho_0) } \subset Q_{\rho_0} \ \ \Longrightarrow \ \  \osc_{Q^{ \omega_0}_{\rho_0,\theta(\omega_0, \rho_0) }} u \leqslant \omega_0  \ . \]

\noindent The H\"older continuity will be a consequence of the iterative procedure leading to the reduction of the oscillation, namely

\begin{propo}\label{pr4.1}
There exist $\epsilon$, $\sigma \in (0, 1)$ such that  the sequences
$\rho_j:=\epsilon^j \rho_0$, $\omega_j:=\sigma^j \omega_0$, $j=0, 1, 2, ...$  satisfy the following assertions 
\begin{equation}\label{eq4.10}
\begin{cases}
Q^{\omega_j}_{\rho_j, \theta_j}:=K^{\omega_j}_{\rho_j}(x_0)\times (t_0- \theta_j , t_0), \ \theta_j =b\,\rho^{p_+}_j \omega^{1-\lambda_+}_j , \quad Q^{\omega_{j+1}}_{\rho_{j+1}}\subset Q^{\omega_j}_{\rho_j},\\
\text{and}\\
\mu^+_j:=\sup\limits_{Q^{\omega_j}_{\rho_j, \theta_j}}u,\quad \mu^-_j:=\inf\limits_{Q^{\omega_j}_{\rho_j, \theta_j}} u,\quad \mu^+_j-\mu^-_j\leqslant \omega_j.
\end{cases}
\end{equation}
\end{propo}

\noindent The proof (by induction) of this result reveals significant differences between the doubly degenerate regime ($\lambda_+>1 \  \wedge \  p_+>2$) and mixed regime  ($\lambda_+>1 \ \wedge \  p_+<2$); we therefore treat them separately.

\begin{remark} 
Observe that the previous construction assures that \eqref{eq4.10} is verified for $j=0$ and also that
\[ \mu_0^+ \leqslant \sup_{Q_{\rho_0}} u \leqslant \omega_0 <(1+\eta) \omega_0, \quad \eta \in (0,1) \ . \] 
We also assume, without loss of generality, that $\mu_0^+ \geqslant \frac{\omega_0}{2}$.
\end{remark}

\subsection{H\"{o}lder Continuity within the doubly degenerate regime $\lambda_+>1 \  \wedge \  p_+>2$}

\noindent Consider the subcylinder
\begin{equation}\label{cylbarQ}
\bar{Q}^{\omega_0}_{\rho_0}:=K^{\omega_0}_{\rho_0}(x_0)\times (t_0- b\rho^{p_+}_0 \omega^{1-\lambda_+}_0, t_0- b\rho^{p_+}_0 \omega^{1-\lambda_+}_0+\,\rho^{p_+}_0\,\omega^{1-\lambda_+}_0]\subset Q^{\omega_0}_{\rho_0, \theta_0} \ , 
\end{equation}
$b\geqslant 1$ to be chosen depending only on the data, for which one of two situations must hold: for some $\nu \in (0, 1)$, to be defined depending only on the data, either 
\begin{equation}\label{eq4.11}
| \bar{Q}^{\omega_0}_{\rho_0}\cap\{u\geqslant \mu^+_0-\frac{1}{4}\,\omega_0\}|\leqslant \nu |\bar{Q}^{\omega_0}_{\rho_0}|,
\end{equation}
or
\begin{equation}\label{eq4.12}
| \bar{Q}^{\omega_0}_{\rho_0}\cap\{u\geqslant \mu^+_0-\frac{1}{4}\,\omega_0\}|\geqslant \nu |\bar{Q}^{\omega_0}_{\rho_0}|.
\end{equation}

\paragraph{Analysis of the First Alternative.}
Assume \eqref{eq4.11} is in force. Then a straightforward application of Lemma \ref{lem2.7} and Remark $2.2$ (taking $b=1$ and $\nu=\gamma^{-1}$) gives
\[ u\leqslant \mu^+_0-\frac{1}{8}\omega_0,\quad \text{a.e.}\quad K^{\omega_0}_{\rho_0/2}(x_0) \times (t_0- \theta_0 + \frac{1}{2}\rho^{p_+}_0 \omega^{1-\lambda_+}_0,  \,t_0- \theta_0 + \rho^{p_+}_0 \omega^{1-\lambda_+}_0 ]  \ . \]
 In particular
\[ u(x,t_0- b\rho^{p_+}_0 \omega^{1-\lambda_+}_0+\,\rho^{p_+}_0\,\omega^{1-\lambda_+}_0) 
\leqslant \mu^{+}_0-\frac{1}{8}\omega_0 \ ,  
\quad \mathrm{in} \quad x \in K^{\omega_0}_{\rho_0/2}(x_0) \ , \]
and then Lemma  \ref{lem2.8} together with Remark $2.3$, for $p_i$'s close to each other, assures the existence of $\xi=\xi(b)\in (0, 1/8)$, such that  
\begin{equation}\label{eq4.13}
\osc\limits_{ Q^{\omega_0}_{\frac{1}{4}\rho_0, \theta(\omega_0, \frac{1}{4}\rho_0)}} u\leqslant \big(1-\frac{1}{4}\xi\big)\,\omega_0.
\end{equation}

\paragraph{Analysis of the Second Alternative.} Assume now that \eqref{eq4.12} holds. Then one also has 
\begin{equation}\label{eq4.14}
| K^{\omega_0/4}_{\rho_0}(x_0)\cap\big\{u(\cdot, \bar{t})\geqslant \frac{1}{4}\omega_0\big\}|\geqslant
| K^{\omega_0}_{\rho_0}(x_0)\cap\{u(\cdot, \bar{t})\geqslant \mu^+_0-\frac{1}{4}\,\omega_0\}|\geqslant \frac{\nu}{4^{N/p}} |K^{\omega_0/4}_{\rho_0}(x_0)|
\end{equation}
for some time level
$t_0-b\,\rho^{p_+}_0\,\omega^{1-\lambda_+}_0 < \bar{t}\leqslant t_0-b\,\rho^{p_+}_0\,\omega^{1-\lambda_+}_0+\,\,\rho^{p_+}_0\,\omega^{1-\lambda_+}_0$. Thereby, and by means of Lemma \ref{lem4.1}, for $C_*, B >1$, depending only on the data,
\begin{equation}\label{eq4.15}
u(x, t)\geqslant \frac{\omega_0}{4C_*},\quad a.e. \ \ (x, t) \in K^{\omega_0/(4C_*)}_{\rho_0}(x_0) \times (\bar{t}+\frac{1}{2}B\,\rho^{p_+}_0\, (\frac{\omega_0}{4})^{1-\lambda_+}, \bar{t}+B\,\rho^{p_+}_0\,(\frac{\omega_0}{4})^{1-\lambda_+}] \  ,
\end{equation}
provided $\lambda_+-\lambda_- \leqslant \epsilon_*$, for some $\epsilon_*= \epsilon_*(data) \in (0,1)$. We then choose $\sigma=\frac{1}{4C_*}$, $\eta= \frac{\sigma}{2}$, and $b=4^{\lambda_+-1} B$ and consider a small positive number $\epsilon$, such that $\epsilon^{p_+} \leqslant \sigma ^{\lambda_+-\lambda_-} $, hence \eqref{eq4.15} yields
\begin{equation}\label{eq4.16}
\osc\limits_{K^{\sigma \omega_0}_{\epsilon \rho_0}(x_0)\times (t_0-b \ \omega_0^{1-\lambda_+} (\frac{\rho_0}{4})^{p_+}, t_0)} u\leqslant (1+\eta-\sigma)\,\omega_0=(1-\frac{1}{2}\sigma)\, \omega_0 \ .
\end{equation}

\paragraph{Reduction of the oscillation.} Combining \eqref{eq4.13} and \eqref{eq4.16} and choosing 
\begin{equation}\label{eq4.17}
\sigma_0= \max \left\{1-\frac{1}{4}\xi,  1-\frac{1}{2} \sigma\right\} > \sigma \quad \mbox{and} \quad
\epsilon^{p_+}= \frac{1}{4^{p_+}}\min\left\{\sigma^{\lambda_+-\lambda_-}, \sigma_0^{\lambda_+-1}\right\}, 
\end{equation}
we obtain
\[
\osc\limits_{Q^{\omega_1}_{\rho_1, \theta(\omega_1, \rho_1)}} u\leqslant \omega_1:=\omega_1,\quad \rho_1:=\epsilon \,\rho_0 \, \ \ \omega_1= \sigma_0\,\omega_0
\]
which proves \eqref{eq4.10} on the first step $j=1$.

\paragraph{The Proof of the H\"{o}lder Continuity Concluded.} We can repeat this procedure, assuming on each step that
\begin{equation}\label{eq4.18}
\mu^+_j\leqslant (1+\frac{1}{2} \sigma)\,\omega_j
\end{equation}
and obtain \eqref{eq4.10} with $\epsilon$, $\sigma_0$ as in \eqref{eq4.17}.  Assume now that $j_0\geqslant 1$ is the first positive integer for which \eqref{eq4.18} is violated, i.e.
\begin{equation}\label{eq4.19}
\mu^+_{j_0}\geqslant (1+\frac{1}{2} \sigma)\,\omega_{j_0},\quad \text{and hence}\quad \mu^-_{j_0}\geqslant \frac{\sigma/2}{1+\sigma/2}\mu^+_{j_0}.
\end{equation}
By the fact that  inequality \eqref{eq4.18} is true for $j=j_0-1$, one has 
\begin{equation}\label{eq4.20}
\sigma/2\,\, \omega_{j_0}\leqslant \frac{\sigma/2}{1+\sigma/2}\mu^+_{j_0}\leqslant \mu^-_{j_0}\leqslant \mu^+_{j_0-1}\leqslant (1+\frac{1}{2} \sigma)\,\omega_{j_0-1}=\frac{1+\sigma/2}{\sigma_0}\omega_{j_0}.
\end{equation}
From here the proof follows exactly the same argument as the one presented in the proof of Theorem \ref{th1.1}: we start by considering new variables and function 
\[\tau=t \omega^{\lambda_+-1}_{j_0},\quad z_i=(x_i-x_{i,0}) \,  \omega_{j_0}^{\frac{\lambda_+-\lambda_i}{p_i}},\quad i=1, ..., N, \quad v(z,\tau) = \frac{u(x,t)}{\omega_{j_0}}\]
and since $v$ in away from zero on the cylinder $Q_{\vec{r}, \eta}(0, \tau_0)$, where $r_i= \rho_{j_0}^{\frac{p_+}{p_i}}$, $\eta= b\rho_{j_0}^{p_+} $ and $\tau_0= t_0 \omega_{j_0}^{\lambda_+-1}$, $v$ satisfies \eqref{eq3.6} in $Q_{\vec{r}, \eta}(0, \tau_0)$, and from here the H\"{o}lder continuity follows.

\subsection{H\"{o}lder Continuity within the mixed regime $\lambda_+>1 \  \wedge \  p_+<2$}

\noindent In order to prove that \eqref{eq4.10} is valid for $j=1$, we study the following alternative: for some $\nu \in (0, 1)$, to be defined depending only on the data, either 
\begin{equation}\label{eq5.12}
| Q^{\omega_0}_{\frac{1}{2}\rho_0, \theta(\omega_0, \rho_0)}\cap\{u\leqslant \frac{1}{4}\,\omega_0\}|\leqslant \nu |Q^{\omega_0}_{\frac{1}{2}\rho_0, \theta(\omega_0, \rho_0)}|
\end{equation}
or
\begin{equation}\label{eq5.13}
|Q^{\omega_0}_{\frac{1}{2}\rho_0, \theta(\omega_0, \rho_0)}\cap\{u\leqslant \frac{1}{4}\,\omega_0\}|\geqslant \nu |Q^{\omega_0}_{\frac{1}{2}\rho_0, \theta(\omega_0, \rho_0)}|.
\end{equation}

\paragraph{Analysis of the First Alternative.} Assuming that \eqref{eq5.12} holds one can conclude that 
\begin{equation*}
u\geqslant \frac{1}{8}\,\omega_0,\quad \text{in}\quad  Q^{\omega_0}_{\frac{1}{4}\rho_0, \frac{1}{2} \theta(\omega_0,\rho_0)},
\end{equation*}
by applying Lemma \ref{lem2.6}, with $\xi=\frac{1}{4}$, $a=\frac{1}{2}$ and $\nu=\frac{1}{\gamma} b^{\frac{N}{p}}$, and making use of Remark $2.1$; thereby getting the reduction of the oscillation 
\begin{equation}\label{eq5.14}
\osc\limits_{Q^{\omega_0}_{\frac{1}{4}\rho_0, \frac{1}{2} \theta(\omega_0,\rho_0)}} u\leqslant \frac{7}{8}\,\omega_0.
\end{equation}

\paragraph{Analysis of the Second Alternative.} Under the second alternative \eqref{eq5.13} there holds
\begin{equation}\label{eq5.15}
| K^{\omega_0}_{\frac{1}{2}\rho_0}(x_0)\cap\{u(\cdot, \bar{t})\leqslant \frac{1}{4}\,\omega_0\}|\geqslant \nu |K^{\omega_0}_{\frac{1}{2}\rho_0}(x_0)| \ ,  \ \ \nu=\frac{1}{\gamma} b^{\frac{N}{p}} \ ,
\end{equation}
for some time level
$t_0-b\,\rho_0^{p_+}\,\omega^{1-\lambda_+}_0\leqslant \bar{t}\leqslant t_0$. From here we want to derive information on the measure of the cubes where $u$ is close to the supremum $\mu_0^+$, for all time levels; for that purpose we use Proposition \ref{pr6.2} in the cylinder
$\bar{Q}^{\omega_0}_{\rho_0}$,
and by \eqref{eq5.15} we obtain
\begin{eqnarray*}
\frac{1}{\gamma}b^{\frac{N}{p}}\,\omega_0& \leqslant &\nu\,(\mu^+_0-\frac{1}{4}\,\omega_0) \leqslant  \fint\limits_{K^{\omega_0}_{\frac{1}{2}\rho_0}(x_0)}(\mu^+_0-u(x, \bar{t}))\,dx \\
& \leqslant\gamma & \inf\limits_{t_0-b\rho^{p_+}_0 \omega^{1-\lambda_+}_0\leqslant t \leqslant t_0}\fint\limits_{K^{\omega_0}_{\rho_0}(x_0)}(\mu^+_0-u(x, t))\,dx
+\gamma\,\omega_0 \sum\limits_{i=1}^N b^{\frac{1}{2-p_i}}\\
& = & \gamma \inf\limits_{t_0-b\rho^{p_+}_0 \omega^{1-\lambda_+}_0\leqslant t \leqslant t_0}\fint\limits_{K^{\omega_0}_{\rho_0}(x_0)}(\mu^+_0-u(x, t))\,dx+
\gamma \omega_0\,b^{\frac{N}{p}+\frac{\bar{\varkappa}}{p(2-p)}}\sum\limits_{i=1}^N b^{\frac{1}{2-p_i}-\frac{1}{2-p}}.
\end{eqnarray*}
Choosing $b$ from the condition
\[ 2N\gamma b^{\frac{\bar{\varkappa}}{p(2-p)}}\leqslant \frac{1}{2 \gamma}\ , \]
and assuming that $\epsilon^{(1)}_*$ satisfies
\[\sum\limits_{i=1}^N b^{\frac{1}{2-p_i}-\frac{1}{2-p}}\leqslant N\,b^{-\frac{p_+-p_-}{(2-p)(2-p_+)}}\leqslant N\,b^{-\frac{\epsilon^{(1)}_*}{(2-p)(2-p_+)}}\leqslant 2N \ ,\]
we obtain, for any $t_0-b \rho^{p_+}_0\,\omega^{1-\lambda_+}_0\leqslant t\leqslant t_0$, 
\[ \frac{1}{2 \gamma} b^{\frac{N}{p}}\omega_0\leqslant \gamma\,\fint\limits_{K^{\omega_0}_{\rho_0}(x_0)}(\mu^+_0-u(x, t))\,dx\leqslant 
\gamma \epsilon \omega_0 + 2 \gamma \omega_0\frac{|K^{\omega_0}_{\rho_0}(x_0) \cap\{\mu^+_0-u(\cdot, t)> \epsilon\,\omega_0\}|}{|K^{\omega_0}_{\rho_0}(x_0)|} \ \]
and hence, choosing $\epsilon$ such that $ \gamma \,\epsilon\leqslant \frac{1}{4\gamma}b^{\frac{N}{p}}$, we get
\begin{equation}\label{eq5.16}
|K^{\omega_0}_{\rho_0}(x_0) \cap\{\mu^+_0- u(\cdot, t)> \epsilon\,\omega_0\}| \geqslant \alpha |K^{\omega_0}_{\rho_0}(x_0)| \ , \qquad \alpha= \frac{1}{8\gamma^2}b^{\frac{N}{p}}
\end{equation}
for all $t_0-b \rho^{p_+}_0\,\omega^{1-\lambda_+}_0\leqslant t\leqslant t_0$. And then one also has 
\begin{equation}\label{eq5.17}
|K^{\omega_0}_{\rho_0}(x_0)\cap\{\mu^+_0- u(\cdot, t_0-b\,e^{-\tau(2-p_+)}\,\rho^{p_+}_0 \omega^{1-\lambda_+}_0)\geqslant  e^{-\tau+\tau_0}\,\epsilon \omega_0\}|\geqslant \alpha|K^{\omega_0}_{\rho_0}(x_0)|,
\end{equation}
for all $\tau\geqslant \tau_0$, where $\tau_0>0$ is to be chosen. 

\noindent Introduce the change of variables and the new unknown function
\[\frac{t-t_0}{b \, \rho^{p_+}_0 \omega^{1-\lambda_+}_0}=:- \,e^{-\tau(2-p_+)} \ , \quad y_i:=\frac{x_i-x_{0, i}}{\rho^{\frac{p_+}{p_i}}_0}\, \omega_0^{\frac{\lambda_+-\lambda_i}{p_i}} \ , \quad v:=\frac{e^\tau}{\omega_0}(\mu^+_0-u) \geqslant 0.\] In this new setting, \eqref{eq5.17} translates into
\begin{equation}\label{eq5.20}
| \mathbb{K}_1\cap\{v(\cdot, \tau)\geqslant \epsilon e^{\tau_0}\}|\geqslant \alpha|\mathbb{K}_1| \ , \quad \forall \tau\geq \tau_0 \ ,
\end{equation}
where $\mathbb{K}_1= \{y: |y_i|<1, \forall i\}$ stands for the standard (isotropic) cube centered at the origin. Observe that, in the cylinder $\mathbb{K}_1 \times(\tau_0, \tau_0+e^{\tau_0(2-p_+)})$, $v$ satisfies
\begin{equation}\label{eq5.18}
v_\tau \geqslant \sum\limits_{i=1}^N D_i\bar{a}_i(y, \tau, v, D v) \ , 
\end{equation}
where the functions $\bar{a}_i$, $i=1,\cdots, N$,  satisfy the following structure conditions 
\begin{eqnarray*}
\sum\limits_{i=1}^N \bar{a}_i(y, \tau, v, D v) D_i v & \geqslant & \bar{K_1} \, \sum\limits_{i=1}^N e^{\tau(p_+-p_i)} \left( \frac{ \mu^+_0-e^{-\tau}\,\omega_0\,v}{\omega_0} \right)^{\lambda_i-p_i+1}|D_i v|^{p_i} \\ 
& \geqslant & \bar{K_1}\sum\limits_{i=1}^N  \left( \frac{\mu^+_0-e^{-\tau}\,\omega_0\,v}{\omega_0}\right)^{\lambda_i-p_i+1} |D_i v|^{p_i},\\
|\bar{a}_i(y, \tau, v, D v)| & \leqslant & \bar{K_2} \, e^{\tau(p_+-p_i)} \left( \frac{\mu^+_0-e^{-\tau}\,\omega_0\,v}{\omega_0}\right)^{\lambda_i-p_i+1} |D_i v|^{p_i-1} \\
& \leqslant & 2\bar{K_2} \left( \frac{\mu^+_0-e^{-\tau}\,\omega_0\,v}{\omega_0}\right)^{\lambda_i-p_i+1} |D_i v|^{p_i-1}
\end{eqnarray*}
provided that 
\[ e^{(\tau_0+e^{\tau_0(2-p_+)})(p_+-p_-)}\leqslant e^{(\tau_0+e^{\tau_0(2-p_+)})\epsilon^{(1)}_*}\leqslant 2 \ \]
and considering $\bar{K_1}= b \, (2-p_+) \, K_1 \min\{1, m_-\}$ and $\bar{K_2}= b \, (2-p_+) \, K_2 \max\{1, m_-\}$. Our next goal is to establish a relation, from below and above, between $\mu^+_0-e^{-\tau}\,\omega_0\,v$ and $\omega_0$ which  can be accomplished by working over sets where $v \leqslant k_0= \epsilon e^{\tau_0}$:
\[\frac{\omega_0}{4} \leqslant \omega_0(\frac{1}{2}-\epsilon) \leqslant \mu_0^+ -\epsilon \omega_0 \leqslant  \mu^+_0-e^{-\tau}\,\omega_0\,v \leqslant  \mu_0^+ \leqslant \omega_0 \ .\]
In this setting, the previous structure conditions are $m_i$ independent, meaning that we have an anisotropic p-Laplacean behaviour, for $p_+<2$ (and all $m_i=1$), and \eqref{eq5.20} allows as to say (taking $\epsilon$ smaller if needed)
\[ |\tilde{K_1}^{k_0} \cap \{v(\cdot, \tau_0) \geqslant k_0\}| \geqslant\frac{\alpha}{2} |\tilde{K_1}^{k_0}|\]
where
\[ \tilde{K_1}^{k_0}= \left\{ y: |y_i| \leqslant \frac{1}{k_0^{\frac{p_+-p_i}{p_i}}}, i=1,\cdots, N \right\} \]
provided $\epsilon_*^{(1)}$ is chosen sufficiently small. We now apply Lemma \ref{lem3.2} to conclude for the existence of $C_*>1$ and $\bar{\epsilon}_*$, $\bar{\delta}_*\in (0, 1)$, depending only on the data and $b$, such that
\begin{equation*}
v(y, \tau)\geqslant \frac{k_0}{C_*},\quad (y, \tau)\in \tilde{K_1}^{k_0} \times [ \tau_0+(1-\bar{\epsilon}_*)\bar{\delta}_*\,\epsilon^{2-p_+} e^{\tau_0(2-p_+)}, \tau_0+\bar{\delta}_*\,\epsilon^{2-p_+}e^{\tau_0(2-p_+)}].
\end{equation*}
Returning to the original variables this inequality implies
\begin{equation}\label{eq5.21}
u(x, t_1)\leqslant \mu^+_0-\frac{\epsilon e^{\tau_0-\tau_1}}{C_*}\,\omega_0=\mu^+_0-\xi_1\,\omega_0,\quad  a.e. \ \  x\in K^{\omega_0}_{\rho_0}(x_0),
\end{equation}
for a time level 
$$t_1=t_0-b\,e^{-\tau_1(2-p_+)}\rho^{p_+}_0 \omega^{1-\lambda_+}_0,\quad \tau_1=\tau_0+\bar{\delta}_*\,\epsilon^{2-p_+}e^{\tau_0(2-p_+)}.$$
We are one step away of getting the reduction of the oscillation. In fact, by means of Lemma \ref{lem2.9}, a DeGiorgi type lemma involving initial data, there exists $\nu_0 \in(0,1)$, depending on the data, such that 
\begin{equation}\label{eq5.22}
u(x, t)\leqslant \mu^+_0-\frac{1}{4}\,\xi_1\omega_0,\quad (x, t)\in K^{\omega_0}_{\frac{1}{4}\rho_0}(x_0)\times(t_1, t_1+\nu_0\rho^{p_+}_0 \xi_1^{2-p_+}\omega^{1-\lambda_+}_0).
\end{equation}
Now we choose $\tau_0$ such that
$$t_1+\nu_0\rho^{p_+}_0 \xi_1^{2-p_+}\omega^{1-\lambda_+}_0=t_0,$$
i.e.
$$e^{\tau_0(2-p_+)}=\frac{b}{\nu_0}\Big(\frac{C_*}{\epsilon}\Big)^{2-p_+}.$$
Then \eqref{eq5.22} yields
\begin{equation*}
u(x, t)\leqslant \mu^+_0-\frac{1}{4}\,\xi_1\,\omega_0,\quad (x, t)\in K^{\omega_0}_{\frac{1}{4}\rho_0}(x_0)\times(t_0-
b\,e^{-\tau_1(2-p_+)}\rho^{p_+}_0\,\omega^{1-\lambda_+}_0, t_0),
\end{equation*}
which implies
\begin{equation}\label{eq5.23}
\osc\limits_{K^{\omega_0}_{\frac{1}{4}\rho_0}(x_0)\times(t_0-
b\,e^{-\tau_1(2-p_+)}\rho^{p_+}_0\,\omega^{1-\lambda_+}_0, t_0)} u\leqslant (1-\frac{1}{4}\xi_1)\,\omega_0.
\end{equation}

\paragraph{Reduction of the oscillation }
Combining \eqref{eq5.14} and \eqref{eq5.23} and choosing 
\begin{equation}\label{eq5.24} \sigma=1-\frac{1}{4}\,\xi_1,\quad\epsilon^{p_+}=\frac{1}{4^{p_+}}\,e^{-\tau_1(2-p_+)}\,\sigma^{\lambda_+-\lambda_-} \min\{ 1, \sigma^{\lambda_--1} \}\ ,
\end{equation}
 we arrive at
\[
\osc\limits_{Q^{\omega_1}_{\rho_1}, \theta_1} u\leqslant \omega_1 \ .\]
This proves \eqref{eq4.10} on the first step $j=1$.

\paragraph{The Proof of Theorem \ref{th1.2} Concluded}

We can repeat this procedure, assuming on each step that
\begin{equation}\label{eq5.25}
\mu^+_j\leqslant \frac{17}{16}\,\omega_j
\end{equation}
and obtain \eqref{eq4.10} with $\epsilon$, $\sigma$ as in \eqref{eq5.24}.  Assume now that $j_0\geqslant 1$ is the first number such that \eqref{eq5.25} is violated, i.e.
\begin{equation}\label{eq5.26}
\mu^+_{j_0}\geqslant \frac{17}{16}\,\omega_{j_0},\quad \text{and hence}\quad \mu^-_{j_0}\geqslant \frac{1}{17}\mu^+_{j_0}
\end{equation}
and 
\begin{equation}\label{eq5.27}
\frac{1}{16}\, \omega_{j_0}\leqslant \frac{1}{17}\mu^+_{j_0}\leqslant \mu^-_{j_0}\leqslant \mu^+_{j_0-1}\leqslant \frac{17}{16}\,\omega_{j_0-1}=\frac{17}{16\,\sigma}\omega_{j_0}.
\end{equation}
We now just need to proceed exactly as we did in {\bf{Section 5.1}}: consider new variables and function 
\[ \tau=t \omega^{\lambda_+-1}_{j_0},\quad z_i=(x_i-x_{i,0})\,\omega_{j_0}^{\frac{\lambda_+-\lambda_i}{p_i}},\quad i=1, ..., N \ , \quad v(z,\tau)= \frac{u(x,t)}{\omega_{j_0}} \]
so that $v$ is away from zero in $K_{\vec{r}}(0) \times (\tau_0 -\eta, \tau_0)$ for $\tau_0 = t_0 \omega_{j_0}^{1-\lambda_+}$, $\eta= b \rho_0^{p_+}$ and $r_i= \rho_0^{\frac{p_+}{p_i}}$, which will lead to the H\"older continuity.

\vskip3.5mm

\noindent {\bf Acknowledgements.}  
S. Ciani acknowledges the support of GNAMPA (INdAM), the PNR funding of the Italian government and the department of Mathematics of the University of Bologna Alma Mater. E. Henriques was financed by Portuguese Funds through FCT - Funda\c c\~ao para a Ci\^encia e a Tecnologia - 
within the Project UID/00013/2025, link at https://doi.org/10.54499/UID/00013/2025. Igor Skrypnik and Mariia Savchenko were supported by the National
Research Foundation of Ukraine, project No. 2025.07/0369 “Qualitative methods of nonlinear
analysis of heterogeneous structures.”

\bigskip

CONTACT INFORMATION

\medskip

{\bf Simone~Ciani} (simone.ciani3@unibo.it)\\
Universit\`a di Bologna Alma Mater, Piazza Porta San Donato 5, Italy

\medskip

{\bf Eurica Henriques} (eurica@utad.pt)\\
Centro de Matem\'atica CMAT; Polo CMAT-UTAD\\
Universidade de Tr\'as-os-Montes e Alto Douro, Vila Real, Portugal
  
\medskip 

{\bf Mariia Savchenko} (shan\textunderscore maria@ukr.net)\\ Institute of Applied Mathematics and Mechanics,
National Academy of Sciences of Ukraine, Gen. Batiouk Str. 19, 84116 Sloviansk, Ukraine

\medskip

{\bf Igor I.~Skrypnik} (ihor.skrypnik@gmail.com)\\Institute of Applied Mathematics and Mechanics,
National Academy of Sciences of Ukraine, Gen. Batiouk Str. 19, 84116 Sloviansk, Ukraine\\


\begin{thebibliography}{99}

\bibitem{BogDuzGiaLiaSch}
V. B\"{o}gelein, F. Duzaar, U. Gianazza, N. Liao, C. Scheven,
H\"{o}lder Continuity of the Gradient of Solutions to Doubly Non-Linear Parabolic Equations, Preprint, https://doi.org/10.48550/arXiv.2305.08539.

\bibitem{BurSkr}
K.O. Buryachenko, I. I. Skrypnik, Riesz potentials and pointwise estimates of solutions to
anisotropic porous medium equation, Nonlinear Anal., 178 (2019), 56-85.

\bibitem{CheDiB}
Y.Z. Chen, E. DiBenedetto, H\"{o}lder estimates of solutions of singular
parabolic equations with measurable coefficients, Arch. Rational Mech.
Anal., 118 (1992),  257-271.

\bibitem{CafEva}
L. A. Caffarelli, C. L. Evans, Continuity of the temperature in the two phase Stefan problem, Arch. Ration. Mech. Anal., 81 (1983),  199-220.

\bibitem{CafFri}
L. A. Caffarelli, A. Friedman, Regularity of the free-boundary of a gas
in an $n$-dimensional porous medium, Indiana Univ. Math. J., 29 (1980),  361-391.


\bibitem{CiaGua}
S. Ciani, U. Guarnotta, Liouville rigidity and time-extrinsic Harnack estimates for an
anisotropic slow diffusion, Proceedings of the American Mathematical Society 151(10)
(2023), 4371-4388.

\bibitem{CiaHen}
S. Ciani, E. Henriques, Harnack-type estimates and extinction in finite time for a class of anisotropic porous medium type equations, Calc. Var. 64, 117 (2025). 

\bibitem{CiaHenSavSkr1}
S. Ciani, E. Henriques, M. O. Savchenko, I. I. Skrypnik,
On  local boundedness of  functions from parabolic De Giorgi classes with doubly nonlinear, nonstandard growth
under exact integrability assumptions, Journal of Differential Equations 462 (2026), 114235.

\bibitem{CiaHenSavSkr2}
S. Ciani, E. Henriques, M. O. Savchenko, I. I. Skrypnik,
Qualitative properties of solutions to anisotropic parabolic equations: Part I - Expansion of positivity, J. Evol. Equ. 26, 84 (2026).

\bibitem{CiaHenSavSkr3}
S. Ciani, E. Henriques, M. O. Savchenko, I. I. Skrypnik,
Qualitative properties of solutions to anisotropic parabolic equations: Part II - Anisotropic Trudinger's equation,  J. London Math. Soc., 114: e70643 (2026).

\bibitem{CiaHenSkr24}
S.Ciani, E. Henriques, I. I. Skrypnik, The impact of intrinsic scaling on the rate of extinction for anisotropic non-Newtonian fast diffusion, Nonlinear Analysis 242, (2024)
, 113497

\bibitem{CiaHenSkr}
S. Ciani, E. Henriques, I. I. Skrypnik, On the continuity of solutions to anisotropic elliptic operators in the limiting case, Bull. London Math. Soc., 57: 1548-1567 (2025)

\bibitem{CiaMosVes}
S. Ciani, S. Mosconi, V. Vespri,  Parabolic Harnack Estimates for anisotropic slow diffusion,
 JAMA 149, (2023), 611-642.
 
 \bibitem{CiaSkrVes}
 S. Ciani, I. I. Skrypnik, V. Vespri,  On the local behavior of local weak solutions to some singular anisotropic elliptic equations,
 Adv. in Nonl. Analysis, 12(1) (2023), 237-265.


\bibitem{CVV} S. Ciani, V. Vespri, M. Vestberg. {\it Boundedness, Ultracontractive Bounds and Optimal Evolution of the Support for Doubly Nonlinear Anisotropic Diffusion.} Annali di Matematica Pura ed Applicata (1923-) (2025), 1-51.

\bibitem{DiB_1}
E. DiBenedetto, On the local behaviour of solutions of degenerate parabolic
equations with measurable coefficients, Ann. Sc. Norm. Sup. Pisa Cl. Sc.
Serie IV, XIII 3 (1986), 487-535.

\bibitem{DiB_2}
E. DiBenedetto, Degenerate Parabolic Equations, Universitext, Springer-Verlag, New York,
1993.

\bibitem{DiBFri}
E. DiBenedetto, A. Friedman, H\"{o}lder estimates for nonlinear degenerate parabolic systems, 
J. Reine Angew. Math., 357 (1985),  1-22.

\bibitem{DiBGiaVes1}
E. DiBenedetto, U. Gianazza, V. Vespri,
 Harnack Estimates for Quasi-Linear Degenerate Parabolic Differential Equation.
Acta Math. 200 (2008), 181-209.

\bibitem{DiBGiaVes2}
E. DiBenedetto, U. Gianazza, V. Vespri, 
Forward, Backward and Elliptic Harnack Inequalities for Non-Negative Solutions
to Certain Singular Parabolic Partial Differential Equations, Ann Scuola Norm Sup Pisa Cl Sci (5), 9  (2010), 385-422.

\bibitem{DiBGiaVes3}
E. DiBenedetto, U. Gianazza, V. Vespri, 
Harnack's inequality for Degenerate and Singular Parabolic Equations. New York:
Springer; 2012.

\bibitem{DiBGiaVes4}
E. DiBenedetto, U. Gianazza, V. Vespri,  Remarks on Local Boundedness and Local H\"{o}lder Continuity of Local Weak Solutions to Anisotropic $p$-Laplacian Type Equations, J Elliptic Parabol Eq. 2 (2016), 157-169.

\bibitem{DMV19} 
F.G. D\"uzg\"un, S. Mosconi, V. Vespri, Anisotropic Sobolev embeddings and the speed of propagation for parabolic equations, Journal of Evolution Equations 19(3) (2019), 845-882.

\bibitem{Hen1}
E. Henriques, Regularity for the porous medium equation with variable exponent: the singular case, Journal of Differential Equations 244 (2008), 2578-2601.

\bibitem{Hen2}
E. Henriques, Concerning the regularity of the anisotropic porous medium equation, Journal
of mathematical analysis and applications 377(2) (2011), 710-731.

\bibitem{Iva1}
A. V. Ivanov, Uniform H\"{o}lder estimates for generalized solutions of quasilinear
parabolic equations that admit double degeneration. Algebra i Analiz, 3(2), (1991), 139-179. Translation in St. Petersburg Math. J. 3 (1992), no. 2, 363-403.

\bibitem{Iva2}
A. V. Ivanov. H\"{o}lder estimates for a natural class of equations of fast diffusion
type, Zap. Nauchn. Sem. S.-Peterburg. Otdel. Mat. Inst. Steklov. (POMI), 229(11), (1995), 299-322. Translation in
J. Math. Sci. (New York) 89 (1998), no. 6, 1607-1630.

\bibitem{Kol1}
I. M. Kolodii, The boundedness of generalized solutions of elliptic differential equations,
Moskow University Vestnik, 25 (1970), 44-52 (In Russian).

\bibitem{Kol2}
I. M. Kolodii,  On boundedness of generalized solution of parabolic differential equations. Moscow University Vestnik, 5,  (1971), 25-31 (In Russian).

\bibitem{LiaSkrVes}
N. Liao, I. I. Skrypnik, V. Vespri,
Local regulay for an anisotropic elliptic equation, Calc.
of Var. and PDE, 59(4), (2020), 1-31.

\bibitem{LisSkr}
V. Liskevich, I. I. Skrypnik, H\"{o}lder continuity of solutions to an anisotropic elliptic equation,
Nonlinear Anal., 71(5-6), (2009), 1699-1708.

\bibitem{NAV26} A. Nastasi, E. P. Ayala, M.  Vestberg, Existence, comparison principle and uniqueness for fully nonlinear anisotropic evolution equations, Nonlinear Analysis, 271, (2026), 114141,

\bibitem{PorVes}
M. Porzio, V. Vespri. H\"{o}lder estimates for local solutions of
some doubly nonlinear degenerate parabolic equations. J. Differential Equations,
103(1), (1993), 146-178.

\bibitem{Ves1}
 V. Vespri, On the local behaviour of solutions of a certain class of doubly
nonlinear parabolic equations. Manuscripta Math., 75(1), (1992), 65-80.

\bibitem{Ves2}
V. Vespri, Harnack type inequalities for solutions of certain doubly nonlinear
parabolic equations. J. Math. Anal. Appl., 181(1), (1994), 104-131.

\bibitem{Tro}
M. Troisi, Teoremi di inclusione per spazi di Sobolev non isotropi, Ricerche Mat. 18 (1969) 3-24.

\end{thebibliography}
\end{document}